%&amstex
\input amsppt.sty
\loadbold
\TagsOnRight
\hsize 30pc
\vsize 47pc
\magnification 1100
\redefine\al{{\overline a}_l}
\redefine\g{\goth g}
\redefine\su{\goth{su}}
\redefine\so{\goth{so}}
\redefine\SU{{\operatorname{SU}}}
\redefine\U{\operatorname{U}}
\redefine\SO{{\operatorname{SO}}}
\redefine\Sp{{\operatorname{Sp}}}
\redefine\Spin{{\operatorname{Spin}}}
\redefine\a{\goth a}
\redefine\B{{\Cal B}}
\redefine\D{{\Cal D}}

\define\Hom{\operatorname{Hom}}

\redefine\l{\goth l}
\redefine\k{\goth k}
\redefine\h{\goth h}
\redefine\m{\goth m}
\redefine\n{\goth n}
\redefine\p{\goth p}
\redefine\t{\goth t}
\redefine\z{\goth z}

\define\q{{\goth q}}
\redefine\C{\Bbb C}
\redefine\R{\Bbb R}

\define\Z{\Bbb Z}

\define\={\overset\operatorname{def}\to=}

\define\CP{\C P}
\define\Mreg{M_{\operatorname{reg}}}

\topmatter
\title Running after a new K\"ahler-Einstein metric \endtitle
\author  Fabio Podest\`a and Andrea Spiro \endauthor 
%\address 
%\endaddress 
%\email 
%\endemail 
\thanks The authors thank the Erwin Schr\"odinger Institut, Wien, for the 
support of their stay in October 1999, when they began working on this project.
\endthanks 
 \abstract We deal with compact K\"ahler manifolds $M$ which are acted on by a
semisimple 
compact Lie group $G$ of isometries with codimension one regular orbits. We
provide an explicit description 
of the standard 
blow-ups of such manifolds along complex singular orbits, in case $b_1(M) = 0$
and the 
regular orbits are Levi nondegenerate. Up to very few exceptions, 
all the nonhomogeneous manifolds in this
class are
shown to admit a $G$-invariant K\"ahler-Einstein metric, giving completely new
examples of
compact K\"ahler-Einstein manifolds.   
\endabstract

\leftheadtext{\smc Running after a new K\"ahler-Einstein metric}
\rightheadtext{\smc F. Podest\`a - A. Spiro}

\endtopmatter
\document

\subhead 1. Introduction
\endsubhead
\bigskip
We consider the class of compact K\"ahler manifolds $M$ with the following two
properties: (a) the first 
Betti number $b_1(M) = 0$; (b) a compact semisimple Lie group $G$ of
(holomorphic) isometries of $M$ acts 
with codimension one regular orbits. Such manifolds, which we will call {\it
K-manifolds\/} throughout the following, 
have been already considered by several authors: many facts on the structure of
$K$-manifolds have been so far discovered and  
successfully used to provide interesting new examples of K\"ahler-Einstein
manifolds (see amongst others \cite{HS}, \cite{Sa}, 
\cite{KS}, \cite{DW}, \cite{PS}, \cite{GC}).\par 
Besides some new general results on K-manifolds, which are contained in first
sections of this paper, we may summarize
our main result in the following Theorem.\par
\bigskip
\proclaim{Main Theorem}ÝLet $M$ be a non-homogeneous, 
K-manifold acted on by 
the semi\-simple Lie group $G$, with only one complex singular $G$-orbit $S$.\par
If $S$ has  
complex codimension one and 
if the regular $G$-orbits are Levi non-degenerate
then: 
\roster 
\item  $M$ is a $G$-homogeneous
holomorphic bundle $G\times_{G_Q, \rho} F$
over the flag manifold $G/G_Q$, 
where $G$, $F$, $G_Q$ and the 
group $Q = \rho(G_Q)$, given by the action $\rho$ of $G_Q$ on the fiber $F$, 
 are as in 
following table:
\medskip
\moveright - 0.1 cm
\vbox{\offinterlineskip
\halign {\strut\vrule\hfil\ $#$\ \hfil
 &\vrule\hfil\ $#$\
\hfil&\vrule\hfil\ $#$\
\hfil&
\vrule\hfil\ $#$\
\hfil&\vrule\hfil\  $#$
\hfil 
\vrule\cr
\noalign{\hrule} n^o &
\phantom{\frac{\frac{1}{1}}{\frac{1}{1}}}G
\ \ &
F
&
G_Q
&
Q
\cr \noalign{\hrule}
1 
&
\underset{\phantom{B}}\to{
\overset{\phantom{A}}\to{
{\SU_n}
}}
&
\underset{\phantom{B}}\to{
\overset{\phantom{A}}\to{
\matrix
Q^2 = \CP^1\times \CP^1\\
 \operatorname{or}\ \CP^2
\endmatrix
}}
&
\operatorname{S}(\operatorname{U}_2 \times \operatorname{U}_{n-2})
&
\SO_3
\cr 
\noalign{\hrule }
2 
&
\underset{\phantom{B}}\to{
\overset{\phantom{A}}\to{\matrix \SU_p\times \SU'_q\\
\smallmatrix p + q > 4\endsmallmatrix
\endmatrix}}
&
\CP^3
&
\operatorname{S}(\operatorname{U}_2\times \operatorname{U}_{p-2})\times 
\operatorname{S}(\operatorname{U}_2\times \operatorname{U}_{q-2}) 
&
\frac{\SO_4}{\Z_2}
\cr \noalign{\hrule}
3
&
\SU_n\ \smallmatrix n >4 \endsmallmatrix
&
\CP^5
&
\underset{\phantom{A}}\to
{\overset{\phantom{B}}\to{
\operatorname{S}(\operatorname{U}_4\times \operatorname{U}_{n-4})
}}
&
\frac{\SO_6}{\Z_2}
\cr \noalign{\hrule}
4
&
\underset{\phantom{A}}\to
{\overset{\phantom{B}}\to{
\SO_{10}
}}
&
Q^7\ \text{or}\ \CP^7
&
\underset{\phantom{A}}\to
{\overset{\phantom{B}}\to{
\SO_2 \times \SO_8
}}
&
\SO_8\ ,\ \frac{\SO_8}{\Z_2}
\cr \noalign{\hrule}
5
&
\underset{\phantom{A}}\to
{\overset{\phantom{B}}\to{
\operatorname{E}_6
}}
&
Q^9\ \text{or}\ \CP^9
&
\underset{\phantom{A}}\to
{\overset{\phantom{B}}\to{
\SO_2\times \operatorname{Spin}_{10}
}}
&
\SO_{10}\ ,\ \frac{\SO_{10}}{\Z_2}
\cr \noalign{\hrule}
}}
\smallskip
\centerline{\bf Table 1}
\medskip
\noindent when $F$ is the complex quadric $Q^r = \SO_{r+2}/\SO_2\times \SO_r$,
the group $Q=\SO_{r+1}$ acts as the standard subgroup of 
$\SO_{r+2}$; when $F$ is the projective space 
$\CP^r =\SU_{r+1}/\operatorname{S}(\operatorname{U}_1\times \operatorname{U}_r)$,
the group $Q=\SO_{r+1}/\operatorname{Center}$ acts as the standard subgroup of
$\operatorname{PSU}_{r+1}$;
\item if $M$ is one of the manifolds described in the Table, 
with the exception of case 1 with $G = \SU_3$
and $F = \CP^2$, case 2 with $G = \SU_p \times \SU_2$, $p>2$, 
case 4 with $F = Q^7$ and case 5 with $F = Q^9$, then 
it is K\"ahler-Einstein with positive first Chern class.
\endroster
\endproclaim
\bigskip
We remark that we can always suppose, up to blow up, that the complex 
singular orbits of a K-manifold are of complex codimension one; moreover the case 
when both singular orbits of a K-manifold are complex can be completely described 
along the lines developed in [PS] (see also [Sp]).\par
At the best of our knowledge all these manifolds are new examples of
nonhomogeneous K\"ahler-Einstein manifolds,
with the only exception of the manifold listed in n.1 with $F = Q^2$, which has
been first discussed by Guan and 
Chen in \cite{GC}. It remains to be checked which of these manifolds 
admits a K\"ahler-Einstein blow down; the present paper already contains 
several results which could be used for such further investigation. \par
The paper is organized as follows. In section 2 we review some basic facts on
K-manifolds and compact 
Levi nondegenerate homogeneous CR-manifolds. In section 3, starting from known
results in \cite{HS} and 
\cite{PS}, we give a fine description of the canonical blow-up of a K-manifold
in Theorem 3.1, while in 
Corollary 3.5 we give the full list of such blow-ups when the regular $G$-orbits
are Levi nondegenerate. \par
In section 4 we describe a generic $G$-invariant K\"ahler metric in terms of a
suitable curve in the Lie 
algebra of $G$ and we write down the Einstein equation for a $G$-invariant
K\"ahler metric for any 
K-manifold. When the regular $G$-orbits are Levi nondegenerate and there exists
only one complex 
singular $G$-orbit of complex codimension one, we also give necessary and
sufficient conditions in order that 
a $G$-invariant K\"ahler-Einstein metric 
on the regular part of $M$ extends as a smooth metric on the whole $M$ (Theorem
4.2).\par
We conclude the paper with section 5, where we prove the existence of a
$G$-invariant K\"ahler-Einstein 
metric on each of the manifolds listed in the main Theorem; this is achieved
proving the existence 
of a solution of the Einstein equation with the appropriate boundary conditions
determined in section 4. 
The proof of this last fact has been inspired by the methods used by Guan and
Chen
 in their paper \cite{GC}. \par
At the moment the authors are not able to see whether the condition of Levi non 
degeneracy is essential for the existence of K\"ahler-Einstein metric on 
K-manifolds with only one complex singular orbit of complex codimension one.
We also stress the fact that our theorem does not state that every  excluded 
case does not admit a K\"ahler-Einstein metric. 
\par
\medskip
As for notation, if $G$ is a Lie group acting isometrically on a 
Riemannian 
manifold $(M, g)$ and $X\in \g$, 
we will adopt the symbol $\hat X$  to denote the corresponding
 Killing vector field on $M$.\par
The  Lie algebra of a Lie
group  
will be always denoted  by the corresponding gothic letter. For a group $G$
and a Lie algebra $\g$, 
$Z(G)$ and $\z(\g)$  denote the center of $G$
and of $\g$, respectively. For any subset $A$ of $G$
or of the Lie algebra $\g$, $C_G(A)$ and $C_\g(A)$ are
the centralizer of $A$ in $G$ and $\g$, respectively. \par
Finally, for any subspace $\n \subset \g$, $\n^\perp$ denotes the 
orthogonal complement of $\n$ in $\g$ w.r.t. the Cartan-Killing form.\par
\remark{Acknowledgment} The authors like to heartily thank Cristina Giannotti for
her 
valuable and long lasting help in the analytic arguments of this paper.\endremark
\bigskip
\bigskip
\subhead 2.  Fundamentals of K-manifolds 
\endsubhead
\bigskip
\subsubhead 2.1 K-manifolds, KO-manifolds and KE-manifolds
\endsubsubhead
\medskip
A {\it K-manifold\/} is a pair formed by a compact K\"ahler manifold $(M,J,g)$ 
and by a compact semisimple Lie group $G$
acting almost effectively and isometrically
(and hence also biholomorphically)
on $M$, such that:
\roster
\item"i)"  $b_1(M)= 0$;
\item"ii)" $G$ acts  by {\it cohomogeneity 
one}, i.e. the regular orbits of 
 the $G$-action are 
of codimension one in $M$.
\endroster
We will also use the notation 
$\omega(\cdot, \cdot)  = g(\cdot, J \cdot)$ for the
K\"ahler fundamental form of
($M,g,J$) and  
 $\rho = r(\cdot, J\cdot)$
for the corresponding the Ricci form.\par
\medskip
For the general properties of  cohomogeneity one manifolds and of  K-manifolds,
see e.g. \cite{AA}, \cite{AA1}, \cite{BR}, \cite{HS}, \cite{PS}. 
At this moment, we  only need to recall the concept of normal geodesic. \par
If $p\in M$ is a  regular point, let us denote by $L = G_p$
the corresponding isotropy subgroup. 
Since $M$ is orientable, every 
regular orbit $G\cdot p$ is orientable. Hence
we may consider a unit normal vector field $\xi$, defined 
on the subset of regular points $\Mreg$, which is orthogonal 
to every  regular orbit. It is known (see \cite{AA1})
that any integral curve of  $\xi$ 
is a geodesic. 
Any geodesic $\gamma:\R\to M$ which is an integral curve of $\xi$ on
$\Mreg$ is called {\it normal geodesic\/} and it crosses every
$G$-orbit orthogonally. \par
\bigskip
The following Proposition will be a basic tool for the sequel.
\par
\bigskip
\proclaim{Proposition 2.1} \cite{PS} Let $(M, J, g)$ be a K-manifold acted on by
the 
compact semisimple Lie group $G$. Let also $p\in \Mreg$ and 
$L = G_p$ the isotropy subgroup at $p$.
Then:
\roster
\item there exists an element $Z$ (determined up to scaling) so that 
$$\R Z \in C_\g(\l)\cap \l^\perp\ , \quad C_\g(\l + \R Z) = \z(\l) + \R Z\ ;\tag
2.2$$
in particular, the connected subgroup $K \subset G$ with subalgebra $\k = \l + \R
Z$
is the isotropy subgroup of a flag manifold $F = G/K$;
\item the dimension of the subspace  $\a = C_\g(\l) \cap \l^\perp $ is either $1$
or $3$; 
in case $\dim_\R \a = 3$, 
then $\a$ is a 3-dimensional subalgebra isomorphic to $\goth{su}_2$
and there exists a Cartan subalgebra 
  $\t^\C \subset \l^\C + \a^\C \subset \g^\C$ so that 
$\a^\C = \C H_\alpha + \C E_\alpha + \C E_{-\alpha}$ for some 
root $\alpha$ of the root system of $(\g^\C; \t^\C)$.
\endroster
\endproclaim
\bigskip
Note that if  
$\dim_\R \a = 1$  (or $\dim_\R \a = 3$) at some regular point $p$, then the same
occurs 
at any other regular point. So we may consider 
the following  definition.\par
\bigskip
\definition{Definition 2.2} Let $(M, J,g)$ be a K-manifold and 
$L = G_p$ the isotropy subgroup of a regular point $p$. We say that 
$M$ is a K-manifold {\it with ordinary action\/} (or shortly, {\it
KO-manifold\/})
if $\dim_\R \a = \dim_\R(C_\g(\l) \cap \l^\perp) = 1$. \par
In all other cases, we say that $M$ is {\it with extra-ordinary action\/} (or,
shortly, 
{\it KE-manifold\/}). 
\enddefinition
\bigskip
\subsubhead 2.2 The structure of K-manifolds. Standard and non-standard
K-manifolds
\endsubsubhead
\par
\medskip
The following is another crucial property of K-manifolds for studying their
structure 
as $G$-manifolds. It is essentially
a corollary of the results in \cite{HS}.\par
\medskip
\proclaim{Proposition 2.3} Let $(M, J, g)$ be a K-manifold acted on by the 
compact
semisimple Lie group $G$. Then it has exactly two singular orbits,
one of which is complex.
\endproclaim
\demo{Proof} It is known that a compact cohomogeneity one manifold
 has either two singular orbits or  no singular orbit
at all. On the other hand,  if there is no singular orbit, it is known that 
the orbit space $\Omega = M/G$ is diffeomorphic
to $S^1$ (see e.g. \cite{AA},  \cite{Br}). But
this   cannot occur,  
because   $b_1(M) = 0$. \par
Consider now an Iwasawa decomposition of $G^\C$, so that we may write 
$G^\C = G\cdot S$, 
where $S$ is a solvable subgroup of $G^\C$. 
We  recall   that, by compactness, 
the 
complexified Lie group $G^\C$ acts on $M$ with one open orbit and that,
by a result of Ahiezer (see
\cite{Ah}), 
$M$ is projective algebraic. From Borel's fixed point
theorem, it follows that there exists a   point $p_o$ which is fixed
by $S$.  Then the $G$-orbit
$G\cdot p_o$ is  complex since it coincides with $G^\C\cdot p_o$ and it
is singular by dimensional reasons.
\qed\enddemo
\bigskip
As we pointed out in the  proof of Proposition 2.3, a  
K-manifold $M$ is acted on by the
complexified Lie group $G^\C$ which has an open orbit. According to the 
definition of Huckleberry and Snow in \cite{HS},
$M$ is called {\it almost homogeneous\/}. \par
In that paper, the authors
give several information on  a  class
of almost homogeneous manifolds,
which contains all the K-manifolds.\par
In the following theorem, we collect some basic features and  the information 
which are immediately implied by the results of Huckleberry and Snow.\par
\medskip
In the following statement and in the rest of the paper, a fixed point for
the $G$-action is  considered as a singular
complex orbit.\par
\bigskip
\proclaim{Theorem 2.4} Let $(M,J,g)$ be a K-manifold acted on 
by the compact semisimple Lie group $G$. Let also $p\in \Mreg$ and 
$\Omega = G^\C\cdot p = G^\C/H$ be the open orbit. \par
{\rm (1)}\ There exists a unique K-manifold $\tilde M$ acted on by $G$, whose 
singular complex  orbits have complex codimension one,  and admitting 
a $G$-equivariant holomorphic map  $\pi: \tilde M \to M$, which is a blow-down
map
along singular complex  orbits.\par
{\rm (2)}\ If $M$ has two singular complex  orbits, then
$\tilde M$ can be $G$-equivariantly and holomorphically fibered onto a
flag manifold $G^\C/P$
$$\tilde \pi: \tilde M \to G^\C/P\ ,$$
where: a) the standard fiber  is $\C P^1$; b) the isotropy of the flag manifold 
$G^\C/P$ is 
the minimal parabolic subgroup  $P \subset G^\C$ which 
 contains $H$; c) the intersection of 
 $\Omega = G^\C\cdot p = G^\C/H$ with the fiber of $\tilde\pi$ is $\C^*$.\par
{\rm (3)}\ If $M$  has exactly one  singular complex orbit, then
$\tilde M$ can be $G$-equivariantly and holomorphically fibered onto a flag
manifold
$$\tilde \pi: \tilde M \to G^\C/P$$
where: 
\roster
\item "a)" the standard fiber is
 $\C P^n$, $Q^n = \{ [z] \in \C P^{n+1}, {}^tz z = 0\}$, $\CP^n \times \CP^n$ ,
 $G_{2,2m}(\C)$ or
$E III = \operatorname{E}_6/\operatorname{Spin}_{10}\times \SO_2$;
\item "b)" the isotropy of the flag manifold
$G^\C/P$
 is  a minimal parabolic subgroup of $G^\C$
containing $H$;
\item "c)" the intersection of
 $\Omega = G^\C\cdot p = G^\C/H$ with the fiber $P/H$ is 
$\C^n$,  $Q^{(n)} = \{ [z] \in \C^{n+1}, {}^tz z = 1\}$,
$\CP^n \times \CP^n \setminus \{[z],[w]|\ {}^tzw = 0\}$,
$\CP^n\setminus
Q^{n-1}$, 
$\Sp_n(\C)/\Sp_{n-1}(\C)$ or $\operatorname{F}_4(\C)/
\operatorname{Spin}_9(\C)$. \par
In each of these cases,  $P/H$
is  the tangent space of a compact rank one symmetric space.
\endroster
\endproclaim
\bigskip
After the results in \cite{HS}, it is convenient to introduce the following 
definitions concerning the different types of K-manifolds.\par
\bigskip
\definition{Definition 2.5}
For any K-manifold $M$, we will call 
the manifold $\tilde M$ defined in  Theorem 2.4 (1) {\it the canonical blow up of
$M$\/}.\par
We will say that a K-manifold $M$ is {\it standard\/} if it has  two
singular complex orbits and {\it non-standard\/} otherwise.\par
A K-manifold $M$  whose canonical blow up $\tilde M$
admits a $G$-equivariant holomorphic 
fibration $\pi:\tilde M\to Q = G^\C/P$ onto a flag manifold  with 
standard fiber  $\CP^1$, will be called {\it projectable\/}. 
\enddefinition
\bigskip
\remark{Remarks 2.6} \par
(i)
If  $\pi: \tilde M \to Q= G^\C/P$ 
is a $G$-equivariant holomorphic fibration 
onto a flag manifold with fiber $\C P^1$, then 
$P$ is a minimal parabolic subgroup containing  $H  \subset G^\C$.\par
\smallskip
(ii) In \cite{PS}, we introduced the concept of {\it K-manifolds with 
projectable
complex structure\/}. It can be proved that such K-manifolds are necessarily
standard 
and hence that {\it being projectable is not the same of having projectable
complex
structure\/}. \par
On the other hand, if  we restrict ourselves to the class
of K-manifolds with ordinary action,  
$M$ is projectable if and only if it has projectable complex structure (see
later).
\endremark
\bigskip
From Definitions 2.2 and 2.5,  the class of K-manifolds is naturally subdivided
into 
four families: {\it standard  KO-manifolds and  standard KE-manifolds\/}, on one
side, and   
{\it non-standard KO-manifolds and non-standard KE-manifolds\/} on the other
side. \par
 A complete description of the standard KO-manifolds 
has been reached in \cite{PS}. An analogous description 
of standard KE-manifolds can be performed following
the same line of arguments used in \cite{PS}.\par
\bigskip
\bigskip
\subsubhead 2.3  The structural decomposition associated with the CR structure 
of a regular orbit
\endsubsubhead
\par
\bigskip
We recall that a {\it CR structure of codimension $r$\/} on  
a manifold $N$ is a pair $(\D, J)$ formed by a distribution $\D\subset TN$ of 
codimension $r$   and a smooth family 
$J$ of complex structures $J_x:
\D_x \to \D_x$ on the spaces  of the distribution. The
CR structure $(\D, J)$ is called {\it integrable\/}
if the distribution $\D^{10} \subset T^\C N$,  given by 
$+i$-eigenspaces $\D^{10}_x \subset \D^\C_x$ of the complex structure $J$
verifies
$$[\D^{10}, \D^{10}] \subset \D^{10}\ .$$
With this definition, we have that any complex structure $J$ on  $N$ 
can be classified as integrable CR structure of codimension zero.\par 
An integrable CR structure $(\D, J)$ of codimension one is called {\it Levi
non-de\-ge\-ne\-ra\-te\/}
if the underlying distribution $\D$ is a {\it contact distribution\/}.
This means that any  local (contact) 1-form
$\theta$, which defines the distribution (i.e. such that $
\operatorname{ker} \theta =
\D$)
is maximally non-degenerate, that is $(d\theta)^n\wedge\theta \neq 0$.
\smallskip
A  smooth map $\phi: N \to N'$ between two CR manifolds $(N, \D, J)$
and $(N', \D', J')$ is  called
{\it CR  map\/} (or {\it holomorphic map\/}) if: 
a) $\phi_*(\D) \subset \D'$; b) for any $x\in N$, 
$\phi_* \circ J_x = J'_{\phi(x)} \circ \phi_*|_{\D_x}$.  
A {\it CR transformation\/} of $(N, \D, J)$ is 
a  diffeomorphism $\phi: N \to N$ which is also a CR map.\par
\medskip
Any submanifold $S$ of a CR manifold
$(N, \D, J)$  is endowed by the  family of subspaces  $\D^S_x\subset T_xS$  and a
family
of complex structures defined by
$$\D^S_x = \{\ v\in (T_xS \cap \D_x)\ :
\ J v\in (T_xS\cap \D_x)\ \}\ \qquad \ \ J_x = J|_{\D^S_x}\ , 
\qquad x\in S\ .$$
If $\D^S = \bigcup_{x\in S}\D^S_x$ is a distribution, we call $(\D^S, J)$ 
{\it induced CR structure\/}. Note that if $S$ is a hypersurface of a complex CR
manifold
$(N,J)$, then $\D^S$ is always a distribution and  $(\D^S, J)$ is an integrable
CR structure 
 of codimension one.\par
\medskip
Let $(G/L, \D, J)$ be a  homogeneous CR manifold of a compact semisimple
Lie group $G$ and assume that  $\D$
is of codimension one. Then $\g$ has the  $\B$-orthogonal decomposition 
$\g = \l + \n$, where $\n$ can be identified with 
the  tangent space $T_{p_o}(G/L)$, $p_o = eL$,  
via the linear isomorphism 
$$\phi : \n \to T_{p_o}(G/L)\ ,\qquad \phi(X) = \hat X|_{p_o}\ .$$
If we denote by $\m$ the subspace
$$\m = \phi^{-1}(\D_{p_o}) \subset \n\ ,$$
 we get the 
following  orthogonal decomposition of $\g$:
$$\g = \l + \n = \l + \R Z_{\D} + \m\ .\tag 2.3$$
where $Z_\D  \in (\l + \m)^\perp$. 
Notice that, since the decomposition 
is $\operatorname{ad}_\l$-invariant, we have that the element  $Z_\D$
is always in $C_\g(\l)\cap \l^\perp$. \par
\smallskip
Using again the identification map $\phi: \n \to T_{p_o}(G/L)$, we may consider
the complex
structure
$$J: \m \to \m\ ,\qquad J  \= \phi^*(J_{p_o})\ .\tag 2.4$$
Note that   $J$ is uniquely determined by the 
 direct sum decomposition 
$$\m^\C = \m^{10} + \m^{01}\ ,\quad \m^{01} 
= \overline{\m^{10}}\ ,\tag 2.5$$
 where $\m^{10}$ and $\m^{01}$ are the 
$J$-eigenspaces  with  eigenvalues $+i$ and $-i$, respectively.\par
\medskip
\definition{Definition 2.7} Let $(N = G/L, \D, J)$ be a compact homogeneous 
CR manifold with an invariant CR structure $(\D, J)$ of codimension one. Then:
\par
\roster
\item"a)" we call {\it the structural decomposition of $\g$ associated with
$\D$\/}
the orthogonal  decomposition  (2.3), with $\m \simeq \D_{eL}$; 
\item"b)" we call {\it the holomorphic\/} (resp. 
{\it anti-holomorphic\/})
{\it subspace associated with $(\D,J)$\/} the subspace
$\m^{10} \subset \m^\C$ (respectively $\m^{01}
= \overline{\m^{10}}$) defined by (2.5).
\endroster
\enddefinition
\bigskip
We also recall that a $G$-invariant CR structure $(\D, J)$ on 
$G/L$ is integrable if and only if the associated holomorphic subspace
$\m^{10} \subset \m^\C$ is so that 
$$\l^\C + \m^{10} \quad \text{is a subalgebra of }\quad \g^\C\ .\tag 2.6$$
We will refer to (2.6) 
as the {\it integrability condition for the holomorphic subspace 
$\m^{10}$\/}.\par
\bigskip
Let us now consider the  regular orbits of a K-manifold. Note that 
if $G/L = G\cdot p_o$ is a regular orbit of $M$ and if 
$(\D, J)$ is the induced CR structure on $G/L$, then  
$(G/L, \D, J)$ is a compact homogeneous CR manifold. Therefore any regular
point $p_o$ determines a structural decomposition for $\g = \l + \R Z_\D(p_o) +
\m(p_o)$, 
$\l = \g_{p_o}$, and a holomorphic subspace $\m^{10}(p_o)$, which are 
those associated with the induced CR structure of $G\cdot p_o = G/L$.\par
\bigskip
\bigskip
\subhead 3. Non-standard  K-manifolds with Levi non-degenerate $G$-orbits
\endsubhead
\bigskip
\subsubhead 3.1 The global structure of a non-standard K-manifold
\endsubsubhead\par
\bigskip
The first main result of this section is the  proof of the following fact.\par
\bigskip
\proclaim{Theorem 3.1} Let $(M,J,g)$ be a non-standard K-manifold
acted on by the compact semisimple Lie group $G$. Then the canonical 
blow-up $\tilde M$ is $G$-diffeomorphic
to a manifold of the form $G\times_{G_Q} F$ where: 
\roster 
\item"a)" $F$ is a $G_Q$-equivariant compactification of the tangent space
$T(G_Q/N)$ 
of a non-trivial compact rank one symmetric space $G_Q/N$; 
\item"b)" $G/G_Q$ is a flag manifold;
\item"c)" if $M$ is a KE-manifold, $G_Q$ is $\SU_2$ and 
$F$ is either $\CP^2$ (= compactification of $T\R P^2 = T(\SO_3/O_2$))
or $\CP^1\times \CP^1$ (= compactification of $TS^2 = T(\SO_3/\SO_2$)); if
$M$ is a KO-manifold, $F$ is  the standard fiber of the $G$-equivariant 
holomorphic bundle $\tilde \pi: \tilde M \to G^\C/P$ given in Theorem 2.4 (3). 
\endroster
\endproclaim
\bigskip
In order to prove this theorem, we first need to make some observations. \par
Let $\tilde M$ be the canonical blow-up of a non-standard K-manifold $M$ and 
let $\hat \pi: \tilde M \to G^\C/P$ the holomorphic projection given in Theorem
2.4 (3). 
It is clear that if the fiber $F$ is not equal to $\CP^1$  (that is {\it $M$ is
non-projectable\/}), 
then the claim of Theorem 3.1 is immediately verified. 
Therefore,  we prove Theorem 3.1 if we can show the following two facts: 
\roster
\item"i)"  a non-standard K-manifold $M$ is projectable only if it is a
KE-manifold; 
\item"ii)" if $M$ 
is a non-standard 
KE-manifold, then it admits a $G$-equivariant holomorphic fibration $\tilde \pi':
\tilde M
\to G^\C/P'$ with fiber equal to $\CP^2$ or $\CP^1\times \CP^1$. 
\endroster
The proof of these two facts  will be the content of the remaining part of this
section, where 
we prove Proposition 3.2 and Theorem 3.3. The content of  Theorem 3.3 is
practically the 
claim  i). The claim ii) is  an immediate consequence of the claims
of Proposition 3.2 and of Theorem 3.3. \par
\bigskip
Let $M$ be a  KE-manifold, $p$ a regular point and $L = G_p$. 
Consider also the connected component  $F$ of the fixed point set
$\operatorname{Fix}(L)$ which passes through $p$. 
Note that $F$ is a closed, 4-dimensional, totally geodesic submanifold, that 
 $T_p F$ is $J$-invariant and that it coincides with the fixed point set 
of the isotropy representation of $L$. \par
In particular, it follows that $F$
is a complex submanifold and hence K\"ahler. Moreover the group 
 $A = \exp(\a) \subset G$, where $\a = C_\g(\l) \cap \l^\perp = \goth{su}_2$,
acts on $F$ 
by cohomogeneity one (see Proposition 2.1 (2)).  If one could check that $b_1(F)
= 0$, 
then we could conclude that $F$ is a K-manifold. \par
In the following Proposition, we show that  this is actually true and we give a
detailed description
of all the possibilities for $F$ and for the action of $A$ on it.\par
\bigskip
\proclaim{Proposition 3.2} Let $M$ be a KE-manifold, $p$ a regular point 
and  $F$ the connected component of 
$\operatorname{Fix}(L)$ through $p$, where $L = G_p$. \par
Then $F$ is a 4-dimensional KE-manifold acted on by the compact semisimple Lie
group 
$A \simeq (N_G(L)/L)^o$. In particular, the action $\rho$ of $A$ on $F$ 
is one of the following:
\roster
\item $\rho(A) = \SU_2$, $F \cong \C P^2 = \SU_3/\operatorname{U}_2$ 
and the action is 
the standard action 
of $\SU_2$ given by the embedding $\SU_2 \subset \SU_3$. 
\item $\rho(A) = \SU_2$, $F$ is a non-trivial $\CP^1$-bundle over 
$\CP^1$ and the group
acts on $F$ by bundle automorphisms. 
\item $\rho(A) = \SO_3$, 
$F \cong \CP^2 = \SU_3/\operatorname{U}_2$ and  the action is
 given by the standard embedding $\SO_3 \subset \SU_3$. 
\item $\rho(A) = \SO_3$, $F \cong \CP^1\times\CP^1$ and 
$A$ acts
 diagonally on $F$. 
\endroster
In cases (1) and (2), $F$ is  standard (actually, in case (1) one of the two
singular
complex orbits 
is a fixed point); in cases (3) and (4)  $F$ is non-standard.
\endproclaim
\demo{Proof} By the previous remarks, $F$ is a 4-dimensional 
K\"ahler manifolds acted on by the 3-dimensional Lie group $A$
with cohomogeneity one and and extra-ordinary action;
note that $A$
is either $\SU_2$ or $\SO_3$.
We have to show that $b_1(F) = 0$. Indeed we recall that the Albanese
map $\alpha: F \to Alb(F)$ is equivariant; moreover $A^\C$, being semisimple,
acts trivially on $Alb(F)$. On the other hand, $A^\C$ has an open orbit
in $F$ and therefore $Alb(F) = \{0\}$, i.e. $b_1(F) = 0$.\par
We now show that only  cases (1), (2), (3) and (4) may occur.\par
If $A$ has a fixed point, then by  the Cone Theorem  (see \cite{HO}), 
$F \cong \C P^2$ and  case (1) occurs. \par
Suppose now that $F$ has no fixed point and  that
the canonical blow-up of $F$ coincides with $F$.\par
If  $F$ is standard,  by Theorem 2.4 (2), 
$F$ is a $\CP^1$-bundle over $\CP^1$ and the group $A = \SU(2)$
acts on $F$ by bundle automorphisms. This bundle cannot be nontrivial,
since otherwise it would be $A$-equivariantly biholomorphic to $\CP^1\times
\CP^1$ 
with the diagonal action of $A = \SU_2$ and at least one
singular orbit would not be complex.\par
Assume now  that  $F$ is non-standard. 
By Theorem 2.4 (3),  there is an $A$-equivariant holomorphic bundle $F \to
A^\C/P$, 
with fiber equal to
$\C P^1$, $\C P^2$ or $Q_2$. If the base $A^\C/P$ is trivial, the 
possibilities for $F$ are either
 $\C P^2$ or $Q_2 \cong \CP^1\times \CP^1$.   But the second case cannot occur
 because $Q_2$ admits an
 $A$-equivariant holomorphic fibration over the 
diagonal $\CP^1$ and $P$ would not be the smallest parabolic subgroup
containing $H$.\par
If the base $A^\C/P$ is not trivial, then 
$F$ must be a $\CP^1$-bundle over $\CP^1$. Any such bundle is $A$-equivariantly 
diffeomorphic to a homogeneous bundle of the form
$$E_k:= \SU_2\times_{\operatorname{T}^1,\rho_k}\CP^1\ ,$$ 
for some $k\in \Bbb N$,
where $\operatorname{T}^1\subset \SU_2$ acts on $\CP^1$ by means of 
the homomorphism 
$$\rho_k: \operatorname{T}^1\to \SU_2\ ,\qquad \rho_k(e^{i\theta}) = 
{\operatorname{diag}}(e^{ik\theta},e^{-ik\theta})\ .$$
 Notice that for $k = 1$, $E_1$ is 
actually $\CP^1\times \CP^1$, where $\SU(2)$ acts diagonally.\par
The proof is concluded if we can show that actually the cases $k >1$  do not 
occur.  
Given a singular point $q \in E_k$, we may suppose that the singular 
isotropy at
$q$ is given by $\operatorname{T}^1$. Notice that the isotropy representation of
$\operatorname{T}^1$ at $q$
decomposes into the sum of the standard isotropy representation of 
$\operatorname{T}^1$
on
$\SU(2)/\operatorname{T}^1$ and the representation $\rho_k$ on 
$\C P^1$. When $k > 1$,
these two representations
are not equivalent and hence the complex structure on $T_pE_k$ has to 
preserve the tangent space to the singular orbit 
$A\cdot q = \SU(2)/\operatorname{T}^1$.
It follows
that when $k>1$ any singular orbit is complex and this contradicts
our hypothesis. 
\qed
\enddemo
\bigskip
 \proclaim{Theorem 3.3} Let $M$ be a non-standard 
K-manifold, $G$  the  compact semi\-simple 
Lie group acting on $M$ and $L = G_p$, $p\in M$, a regular isotropy subgroup.
Then:
\roster 
\item"a)" $M$ is projectable only if the action of $G$ is extra-ordinary;
\item"b)" assume that $M$ is projectable; if  
$F$ is the connected component of $\operatorname{Fix}(L)$ through $p$,
then 
$F \cong \CP^1\times \CP^1$; furthermore, $\tilde M$ fibers
holomorphically 
and $G$-equivariantly onto a $G^\C/\tilde P$, with 
$\tilde P$ parabolic in $G^\C$ and  
 with standard fiber equal to $F = \CP^1\times \CP^1$. 
\endroster
\endproclaim
\demo{Proof} a) We will show that if $M$ is projectable  
with ordinary action,
then it is 
standard. With no loss of generality, we will assume  that 
$M = \tilde M$, where 
$\tilde M$ is the canonical blow-up.\par
Let $\gamma$ be a normal geodesic through $p$ and $\pi: M \to G^\C/P$
a $G$-equivariant holomorphic fibration onto a flag manifold, 
with fiber $\CP^1$.
Let also  $K = P \cap G$, so that we may write
 $G^\C/P = G/K$. \par
For any regular point $\gamma_t$, the fibration
$\pi$ induces a  CRF map   $\hat \pi: G\cdot \gamma_t = G/L \to G/K$ 
 with fiber $K/L = S^1$. In particular, we have that
$\k = \l + \R Z$
 for some $Z \in C_\g(\l) \cap \l^\perp$. \par
On the other hand, by Prop. 2.1 in \cite{PS},
 the moment map $\mu: M \to \g^*$ induces 
a $G$-equivariant map from  $G/L = G\cdot \gamma_t$ onto a 
flag manifold $G/K_t$
$$\mu_t : G/L = G\cdot \gamma_t \to  G/K_t$$
where $L$ is  of codimension one in $K_t$. In particular
$$\k_t = \l + \R Z_t$$
for some $Z \in C_\g(\l) \cap \l^\perp$. 
Since $C_\g(\l) \cap \l^\perp$ is 1-dimensional we conclude that $\k_t = \k$ 
and that $G/K_t = G/K$ for any $t$. \par
Furthermore, if 
$\g = \l + \R Z_\D(t) + \m_t$ is the structural decomposition 
associated with the induced CR structure of $G\cdot \gamma_t = G/L$, using once
again 
the fact that $\dim_\R C_\g(\l) \cap \l^\perp = 1$, we conclude that
$\k = \k_t = \l + \R Z_\D(t)$, for any regular point $\gamma_t$.\par
To conclude, consider the anti-holomorphic subspace $\m^{01}_t$ of 
the orbit $G\cdot \gamma_t = G/L$. By Proposition 2.8 a) and b), we have that
$$\p = \k^\C + \n = (\l^\C + \C Z_\D(t)) + \n\ ,
\quad \l^\C + \m^{01}_t \subset
\p$$
where $\n$ is the nilradical of $\p$. 
Since $\dim_\C \m^{01} = \dim_\C \n$ and they are both $\B$-orthogonal to 
$\l^\C + \C Z_\D$, it follows that $\m^{01} = \n$. But $\n$ is the nilradical 
and hence 
$$[\k, \m^{01}]  \subset \m^{01} \ ,\qquad [\k, \m^{10}]  =
 \overline{[\k, \m^{01}]} \subset \m^{10}\ .$$
This means that the induced
CR structure $(\D, J)$ of $G\cdot \gamma_t$ is
$\operatorname{ad}_{\k}$-invariant, that is it is {\it projectable\/} in the
sense of 
\cite{PS}. The conclusion follows from the fact that in this case, by Prop. 4.1
in \cite{PS}, 
both singular orbits are complex.\par
\bigskip
b) First of all we prove that $F$ has only one complex $A$-orbit.
 We recall
that $F$ contains a normal geodesic and that a singular $A$-orbit in $F$
is contained in a singular $G$-orbit in $M$. We will show that a complex
$A$-orbit in $F$ is contained in a complex $G$-orbit and 
therefore $F$ has
only one complex orbit since $M$ is non-standard. Let $p\in F$ such that
$A\cdot p$ is a complex orbit; we denote by $N'$ and $N$ the 
normal spaces
to the $A$- and $G$-orbits respectively and by $v\in N\subset N'$ the tangent
vector
of a normal geodesic through $p$. Now let $w$ be
a unit vector in $N$; since the isotropy $G_p$ acts 
transitively on the unit
sphere of $N$, we may find $g\in G_p$ with $gw = v$ and since the normal
space $N'$ is complex, we may find $h \in A_p$ such that $hv = Jv$. Then
$Jw = gJv = ghv$, meaning that $Jw \in N$ and therefore $N$ is complex.\par
It the follows form by Proposition 3.1 that $F$ has no fixed point and it is
either $F \cong \CP^1\times\CP^1$ or $F \cong \CP^2$.  This implies  that 
the canonical blow up of $F$ coincides with $F$ itself. Hence, 
 we may again assume  that $M = \tilde M$ with no loss of generality. \par
Since $M$ is projectable, there exists a $G$-equivariant holomorphic fibration
 $\pi: M \to G^\C/P = G/Q$ with fiber $\CP^1$ onto the flag manifold $G^\C/P =
G/Q$, where $Q = G\cap P$; for any $x\in M$ we denote by $Z_x$ the fiber
$\pi^{-1}(\pi(x))$. Without loss of generality, we may assume that
$L \subset Q$ and since $\dim Q/L = 1$, we may write the Lie algebra
$\q$ of $Q$ as $\q = \l + \R\cdot a$ for some $a \in \a$; this means that
${\hat a}_p \in T_pF \cap T_pZ_p$ and since both $F$ and the the fiber
$Z_p$ are complex, $T_pZ_p \subset T_pF$. This argument applies actually
to any point $y$ along the normal geodesic through $p$, which is
contained in $F$; by the $G$-equivariance of $\pi$ we have that
$T_xZ_x \subset T_xF$ for all regular point $x \in F$, hence for all $x\in F$.
This means that $Z_x \subset F$ for all $x\in F$. Therefore
$F \simeq \CP^1 \times \CP^1$, since $\CP^2$ is not the total space of a
$\CP^1$-fibration.\par
Let us now consider the holomorphic subspace $\m^{10}$  associated with the 
induced CR structure of $G/L = G\cdot p$ and put
$$\a^{10} = \a^\C\cap \m^{10}\ ,\qquad \a^{01} = \overline{\a^{10}} = \a^\C \cap
\m^{01}\ .$$
Let also $\m'{}^{10} = \m^{10} \cap (\a^{10})^\perp$.\par
Notice that the complex isotropy subalgebras $\h = \g^\C_p$ and $\h^\a = \a^\C_p$
of the actions of $G^\C$
and $A^\C$, respectively, are
 equal to 
$$\h = \l^\C + \m^{01} = \l^\C + (\a^{01} + \m'{}^{01})\ ,\quad
\h^\a = \a^{01}\ .$$
Note that  
$\p = Lie(P)$ is a minimal parabolic subalgebra of $\g^\C$ properly containing
$\h$. Since $F$ contains every fiber $Z_x$ of $\pi$ for $x \in F$, we have
that $\dim_\C \p \cap \a^\C = 2$; moreover since
$\a^{01} \subset \p \cap \a^\C$ is generated by a regular element of $\a^\C$,
as it can be checked directly 
using the explicit action of $A$ on $F = \CP^1 \times
\CP^1$, we conclude that $\p \cap \a^\C$ is a parabolic subalgebra of $\a^\C$.
\par
Since
$\dim_\C \p = \dim_\C \h + 1$, it follows that $\p$ is equal to 
$$\p = (\p \cap \a^\C) + \m'{}^{01}\ .$$In particular $\m'{}^{01}\subset \p$. 
We now claim that the subspace 
$$\tilde\p =  \p + \a^\C = \l^\C + \a^\C + \m'{}^{01}$$
is  also a  parabolic subalgebra.\par
First of all, consider a Cartan subalgebra $\t_{\l}^\C$
of  $\l^\C$. By Proposition 2.1, for any  regular element 
$Z\in \a^\C$, the subalgebra 
$\t_{\l}^\C + \R Z$ is a Cartan 
subalgebra for $\g^\C$. We
already remarked that $\a^{01}$ is generated by a regular
element of $\a^\C = \goth{sl}_2(\C)$ and therefore 
$\t^\C = \t_{\l}^\C + \a^{01} \subset \p$ is a Cartan subalgebra for $\g^\C$ 
included in $\p$. \par
Now,  consider 
the root system $R$ of $\g^\C$ determined by $\t^\C$ and  denote by $S$ and $P$
the closed subsystems of $R$ defined by
$$S = \{ \beta\in R\ : \ E_\beta \in \h\ \}\ ,\qquad
P = \{ \beta\in R\ : \ E_\beta \in \p\ \}\ .$$
We recall that 
$\a^\C = \a^{01} + \C\cdot E_\alpha + \C\cdot E_{-\alpha}$ for some
root $\alpha \in R$. Since $\a^\C\cap\h = \a^{01}$, 
it follows that  $\pm\alpha\notin S$. On the other 
hand, since  $P$ is a parabolic subsystem, for any root $\gamma \in R$, either
$\gamma$ or $-\gamma$ is in $P$. Hence we may  assume
$P = S\cup \{\alpha\}$. \par
To prove that $\tilde \p = \p + \a^\C$ is a parabolic subalgebra, we have 
only to check that  the subset of  $\tilde P = S\cup \{\alpha,-\alpha\}$ is
closed. This reduces to show that $(S + \{-\alpha\}) \cap R \subseteq \tilde P$;
suppose not, then 
there exists a root $\beta\in S$ so that $\beta -\alpha\in R\setminus \tilde P$.
In particular 
$\beta - \alpha\notin P$ and since $P$ is parabolic, this implies that $\alpha -
\beta \in P$, 
actually $\alpha - \beta\in S$. But then $\alpha = (\alpha - \beta) + \beta \in
(S+S)\cap R\subset S$
and this is a contradiction. \par
If $\tilde P$ denotes the parabolic subgroup of $G^\C$ with Lie algebra
$\tilde\p$,
then $\tilde M$ 
fibers holomorphically and $G$-equivariantly onto 
$N = G^\C/{\tilde P}$, with a complex
$2$-dimensional 
fiber $\Cal F$. Since $A\subset \tilde P \cap G$, we get that
 $A$ acts almost effectively on $\Cal F$ and hence that
 $F\cap \Cal F$ is at least three dimensional. Since
$\Cal F$ and $F$ are both complex, we conclude that $\Cal F = F$. \qed
\enddemo
\bigskip
\bigskip
\subsubhead 3.2 Non-standard K-manifold with Levi non de\-ge\-ne\-ra\-te orbits 
\endsubsubhead
\par
\bigskip
We  now reduce 
to consider only K-manifold with regular $G$-orbits, which
are  Levi non-degenerate. Notice that if 
one regular $G$-orbit is Levi non-degenerate, then 
all regular $G$-orbits are Levi non-degenerate. \par
The complete list of all simply connected, compact homogeneous manifolds
$G/L$ of a compact semisimple Lie group, which admit a $G$-invariant
Levi non-degenerate integrable CR structure $(\D, J)$ of codimension one,  has
been obtained 
 in \cite{AS} (see also \cite{AHR}). According to the results in \cite{AS}, 
any such simply connected homogeneous CR manifold falls in one of the following
three families:
\roster
\item"a)" $(G/L, \D, J)$ is a homogeneous $S^1$-bundles $\pi: G/L \to F = G/K$
over a flag manifold $F = G/K$ with invariant complex structure
$J_F$,
and  $(\D, J)$ is the unique CR structure such that the map $\pi$ is holomorphic;
\par
\item"b)" $(G/L, \D, J)$ is a 
sphere bundles $G/L = S(N)\subset TN$ of a compact
rank one symmetric space $N = G/H$, with the CR structure $(\D, J)$ 
is induced by
the 
natural complex structure of $TN = G^\C/H^\C$; \par
\item"c)" $G/L$ is one of the following manifolds: $\SU_n/\operatorname{T}^1\cdot
\SU_{n-2}$, 
$\SU_p\times \SU_q/\operatorname{T}^1 \cdot$ $\U_{p-2}\cdot \U_{q-2}$, 
$\SU_n/\operatorname{T}^1\cdot \SU_2\cdot \SU_2\cdot \SU_{n-4}$, 
$\SO_{10}/\operatorname{T}^1\cdot \SO_6$, 
$\operatorname{E}_6/\operatorname{T}^1\cdot \SO_8$;
these manifolds admit canonical holomorphic fibrations over 
a flag manifold $(F,$ $J_F)$ 
with typical fiber $S(S^k)$, where $k = 2, 3, 5, 7, 9$ or $11$, respectively;
the CR structure  is determined by the invariant complex
structure  $J_F$ on $F$ and by 
an invariant CR structure on the typical fiber,
depending on one complex parameter.
\endroster
\bigskip
From Theorem 2.4, Theorem 3.1
 and the above quoted classification of 
compact homogeneous CR manifolds, the following Corollary 
is obtained.\par
\bigskip
\proclaim{Corollary 3.5} Let $(M, J, g)$ be a 
non-standard K-manifold with one regular 
$G$-orbit $G/L = G\cdot x$, which is Levi non-degenerate.
Then only one of the following cases may occur (in what follows, $\tilde M$
is the  blow up of 
$M$ along the unique singular $G^\C$-orbit):  
\roster 
\item"i)" 
$M = \tilde M$ and it
is $G$-equivariantly biholomorphic to one of the following compactifications
of  the tangent space   of a compact rank one symmetric space $S$:\par
a) $\CP^n$; in this case $S = \R P^n$ and 
$G = \SO_{n+1}$ or $G = \operatorname{Spin}_7$
if $n = 7$;\par
b)  $Q^n = \{ [z] \in \C P^{n+1}, {}^tz z = 0\}$; in this case $S = S^n$ and 
 $G = \SO_{n+1}$ or $G = \operatorname{Spin}_7$
if $n = 7$;\par
c) $\CP^n\times \CP^n$; in this case $S = \C P^n$ and $G = \SU_{n+1}$; \par
d) $Gr_{2,2m}(\C)$; in this case $S = \Bbb H P^n$ and $G = \Sp_n$;\par
e) $EIII = \operatorname{E}_6/(\SO_2\times \operatorname{Spin}_{10})$; 
in this case
$S = \Bbb O P^2$ and $G = \operatorname{F}_4$; 
\item"ii)" $\tilde M$ is biholomorphically $G$-equivalent  to a manifold 
of the form $G \times_{G_Q, \rho} F$, where $G$, $F$, $G_Q$ and the 
group $Q = \rho(G_Q)$, given by the action $\rho$ of $G_Q$ on $F$, 
are as in one of the
cases of Table 1. The  cases  in $n^o.$1 of Table 1  are the only 
possibilities for non-standard KE-manifolds; all other cases correspond 
to non-standard KO-manifolds.\endroster
\endproclaim
\demo{Proof} The content of the corollary follows directly 
from the above remarks. 
In particular, Table 1 has been obtained as follows: 
the groups $G$ and $G_Q$ are determined by 
the list given at point c) at the beginning of this 
section; the groups $Q$ are determined as the only groups, 
which are homomorphic images of $G_Q$ and acting non-standardly 
on one of the manifolds listed in (3.a) of Theorem 2.4, and the fibers
$F$ are determined accordingly.\qed
\enddemo 
\bigskip
\remark{Remark 3.6}Ý
We  observe that  each manifold $M = G \times_{G_Q, \rho} F$, where 
$G$, $G_Q$, $Q = \rho(G_Q)$ and $F$ are as in Table 1,
and for each 
of the two $G$-invariant complex structure $J_o$ on the flag manifold $G/G_Q$,
there exists a  $G$-invariant 
complex structure $J$ on $M$, such that the canonical  
projection $\pi: (M, J) \to (G/G_Q, J_o)$
is holomorphic. \par
We indicate how to prove this claim  in case 1 of Table 1.
 The flag manifold $G/G_Q$ can be written as
$G^\C/P = \operatorname{SL}_n(\C)/P$, where 
$P$ is a parabolic subgroup. $P$  admits a holomorphic projection
onto $S = \operatorname{SL}_2(\C)/Z$, $Z$ being the center 
and  $S$ acts holomorphically 
onto $Q^2$ and onto $\C P^2$ in a standard way, with an action $\tilde \rho$,
which extends the  action $\rho$ of $G_Q$
on $F$. Therefore, the manifold $M$ is $G$-equivariantly 
diffeomorphic to $\operatorname{SL}_n(\C) \times_{P, \tilde \rho} F$;  this last
manifold can be shown to be $G^\C$-homogeneous, holomorphic bundle over 
$G^\C/P = G/G_Q$, proving our claim.\par
The cases 2, 3 and 4 can be checked similarly. In case 5, it is enough to check
that the parabolic subgroup $P$ such that $\operatorname{E}_6(\C)/P = 
\operatorname{E}_6/\SO_2\times \Spin_{10}$ admits a projection
onto $\Spin_{10}(\C)$; this can be achieved considering that the isotropy action
of
$P$ on the tangent space of the symmetric space 
$\operatorname{E}_6/\SO_2\times \Spin_{10}$ coincides with the action 
of $\C^* \times \Spin_{10}(\C)$ (see Table 2 in \cite{Be}, p. 313).\par
Notice also that  the claim implies that 
each of the fiber bundles described in Table 1 
does correspond to a K-manifold, 
since each of them can be realized as an algebraic variety.
\endremark
\bigskip
\remark{Remark 3.7} The manifold 
$M = G \times_{G_Q, \rho} F$  described in case 1 of Table 1, with $F = Q^2$,
is the  manifold discussed in \cite{GC}.
\endremark
\bigskip
\bigskip
\subhead 4. The Einstein equation for a non-standard K-manifold
\endsubhead
\bigskip
\subsubhead 4.1 Optimal transversal curves and the algebraic representatives of 
closed 2-forms
\endsubsubhead
\par
\medskip
By the results of \cite{Sp}, it is known that on any K-manifold $M$ of real
dimension $2n$, 
there exists 
a family of curves $\eta: \R \to M$ which verify the following properties:
\roster
\item the points $\eta_t$ are  of the form 
$$\eta_t = \exp(i t Z) \cdot p_o$$
for some  $p_o \in M$, which is regular for the $G^\C$-action and 
for some $Z\in \g$; more precisely, in case $M$ is standard, $p_o = \eta_0$
is any $G$-regular point; in case $M$ is non-standard, $p_o$ is a
point of the singular $G$-orbit, which is not complex; 
\item $\eta$ intersects any regular $G$-orbit; in particular, in case $M$ is
standard, 
then $\eta_t \in \Mreg$ for any $t$; in case $M$ is non-standard, 
$\eta_t \in \Mreg$ if and only if $t \neq 0$;
\item for any  point $\eta_t \in \Mreg$,
 the tangent vector
$\eta'_t = J \hat Z_{\eta_t}$
is transversal to the regular orbit $G\cdot \eta_t$; 
\item
any element $g\in G$ which is in a stabilizer $G_{\eta_t}$, $\eta_t\in
\Mreg$,
fixes pointwise the whole curve $\eta$; in particular, all regular
orbits $G\cdot \eta_t$  are $G$-equivalent to the same
homogeneous space $G/L$;
\item  the structural 
decompositions 
$$\g = \l + \R Z_\D(t) + \m(t)$$
associated with the CR structure of the regular orbits $G/L = G\cdot \eta_t$ 
do not depend on $t$; furthermore, $Z_\D(t) = Z$ for any $\eta_t \in
\Mreg$;
\item there exists a basis $\{F_1, G_1, \dots, F_{n-1}, G_{n-1}\}$ for $\m$ such
that for any 
$\eta_t \in \Mreg$
the complex structure $J_t: \m \to \m$,  induced  by the complex structure 
of $T_{\eta_t}M$,  is of the following form:
$$J_t F_{j} =  \lambda_j(t) G_{j}\ ,\qquad J_t G_{j} =  -
\frac{1}{\lambda_j(t)}F_{j}\ ;
\tag 4.1$$
where the function $\lambda_j(t)$ 
is either one of the functions
  $- \tanh(t)$,  $- \tanh(2t)$, $-\coth(t)$ and 
 $- \coth(2t)$ or it is identically equal to  $1$.
\endroster
Any curve $\eta_t$ which verifies (1) - (5) is called {\it optimal transversal
curve\/}. \par
\bigskip
Consider now  a closed $G$-invariant, $J$-invariant
2-form, which is bounded and defined on the regular points subset
$\Mreg$.\par
If  $\eta: \R \to M$ is  an optimal
transversal curve,  since $\g$ is semisimple and $\varpi$ is $G$-invariant,  
then for any $t\in \R$ there exists a unique 
$\operatorname{ad}_\l$-invariant
endomorphism  $F_{\varpi, t}\in \Hom(\g,\g)$ such that:
$$\B(F_{\varpi, t}(X), Y) = \varpi_{\eta_t}(\hat X, \hat Y)\ ,\qquad X, Y \in \g\
.\tag4.2$$
Using the fact that $\varpi$ is $G$-invariant and closed, 
it is not difficult to realize that for any $X, Y, W\in \g$ 
the following holds
$$F_{\varpi,t}([X,Y]), W)  = [F_{\varpi, t}(X), Y] + [X, F_{\varpi, t}(Y)]\ .$$
This means that $F_{\varpi, t}$ is a derivation of $\g$ and hence of the form
$$F_{\varpi, t} = \operatorname{ad}(Z_\varpi(t))\tag 4.3$$
for some $Z_\varpi(t) \in \a = C_\g(\l) \cap \l^\perp$.\par
We call the curve 
$$Z_\varpi: \R \to C_\g(\l) = \z(\l) + \a\ ,\qquad t \mapsto Z_\varpi(t)\tag
4.3$$
the {\it algebraic representative of the 2-form $\varpi$ along the 
optimal transversal curve $\eta$\/}. \par
Note that the 2-form $\varpi$ can be completely
recovered by its algebraic representative $Z_\varpi(t)$. In 
fact, the following proposition holds.
\bigskip
\proclaim{Proposition 4.1}(\cite{PS},\cite{Sp}) Let $(M, J, g)$ 
be a K-manifold 
 acted on by the compact semisimple Lie group $G$. Let also $\eta_t = \exp(t i
Z_\D)\cdot p_o$ 
be an optimal transversal curve
and  $Z_\varpi : \R \to \z(\l) + \a$
the algebraic representative of a bounded, $G$-invariant, $J$-invariant
closed 2-form $\varpi$ along $\eta$.  Then:\par
\roster
\item if $M$ is  a standard K-manifold or a 
non-standard KO-manifold (equivalently, if either
 $\a = \R Z_\D$, or $\a = \goth{su_2}$ and $M$ is standard), 
then there exists an element $I_\varpi \in \z(\l)$ and 
a smooth function ${f}_\varpi: \R \to \R$ so that
$$Z_\varpi(t) = f_\varpi(t) Z_\D +  I_\varpi\ ;\tag 4.4$$
\item if $M$ is non-standard KE-manifold
(equivalently, if $\a = \goth{su_2}$ and $M$ is non-standard), then there exists 
a Cartan subalgebra $\t^\C \subset \l^\C + \a^\C$ and a root $\alpha$
of the corresponding root system, such  that  $Z_\D \in \R(iH_\alpha)$ 
and $\a = \R Z_\D + \R F_\alpha + \R G_\alpha$, where $F_\alpha$ and 
$G_\alpha$ are as in (4.9') below; furthermore
there exists an element $I_\varpi \in \z(\l)$, a real number $C_\varpi$ and 
a smooth function $f_\varpi: \R \to \R$ so that
$$Z_\varpi(t) = f_\varpi(t) Z_\D + \frac{C_\varpi}{\cosh(t)}
G_\alpha +  I_\varpi\ .\tag 4.4'$$
\endroster
Conversely,  if $Z_\varpi: \R \to C_\g(\l)$ is a curve in $C_\g(\l)$ of the form
(4.4)
 or  (4.4'), 
then there exists a 
unique closed $J$-invariant, $G$-invariant 2-form 
$\varpi$  on $\Mreg$, having $Z_\varpi(t)$ as  algebraic representative;
such 2-form is the 
unique $J$- and  $G$-invariant form which  verifies the following 
 at all points of $\eta_t \in
\Mreg$ and for any $V, W \in \m$
$$\varpi_{\eta_t}(\hat V, \hat W) = \B(Z_\varpi(t), [V, W])\ ,
\qquad \varpi_{\eta_t}(J \hat Z_\D, \hat Z_\D)
= - f'_\varpi(t) \B(Z_\D,Z_\D)\ .\tag4.5$$
\endproclaim
\bigskip
We say that the 2-form $\varpi$ 
 is {\it tame\/} if there exists an optimal transversal curve $\eta$, 
along which the algebraic representative $Z_\varpi(t)$ is of the form (4.4). 
By the previous proposition, in case $M$ is a non-standard KE-manifold, 
$\varpi$ is  tame in case if and only if the constant $C_\varpi$, which 
should appear in the expression (4.4') of $Z_\varpi(t)$, is equal to $0$. 
For any other kind of K-manifold, a  bounded $G$-invariant, $J$-invariant closed 
2-form $\varpi$ is always tame.\par
\medskip
Due to Proposition 4.1, a $G$-invariant K\"ahler metric $g$ on a K-manifold
$M$ is  Einstein, with Einstein constant $c>0$, 
if and only if the algebraic representative $Z_\omega(t)$, $Z_\rho(t)$
of the K\"ahler form $\omega$ and the Ricci form $\rho$, along 
an optimal transversal curve $\eta_t = \exp(i t Z_\D) \cdot p_o$, verify 
$Z_\rho(t) = c Z_\omega(t)$ for any $t \in \R$, which is equivalent to 
$$f_\rho(t) = c f_\omega(t) \ ,\qquad  I_\rho = c I_\omega\ ,\qquad C_\rho =
c C_\omega\ ,\tag 4.6$$
the last equation appearing only if $M$ is a non-standard KE-manifold.  
We call (4.6) {\it the Einstein equations of a K-manifold\/}.  \par
\bigskip
\bigskip
\subsubhead 4.2 The Einstein equations of a non-standard K-manifold
\endsubsubhead
\par
\bigskip
We now want  to write down the Einstein equations for a $G$-invariant 
K\"ahler metric $g$ on a non-standard K-manifold $M$. Since we are interested 
in the manifolds listed in Corollary 3.5, we will always assume that 
the canonical blow up $\tilde M$ 
is of the 
form $\tilde M = G\times_{G_Q, \rho} F$ for a subgroup $G_Q$ and an action $\rho$
of  $G_Q$ on $F$.\par 
In what follows we will always denote by 
$N_F$ 
$$N_F = \frac{1}{2} \dim_\R\m_F = 
\frac{1}{2}(\dim_\R (\rho(G_Q) \cdot x) - 1)\ ,$$
for any regular point $x\in F_{\operatorname{reg}}$ in the standard fiber.
\bigskip
Let $\eta_t = \exp(i t Z_\D) \cdot p_o$ be an optimal transversal curve and 
$\g = \l + \R Z_\D + \m$ the structural decomposition given by the CR structure 
of all regular orbits $G\cdot \eta_t = G/L$.
 By Th. 3.2 in \cite{Sp}, we 
may consider  the following orthogonal decomposition
$$\g = (\l_o + \l_F) + \R Z_\D + (\m_F + \m')\ ,
\qquad \l = \l_o + \l_F\ ,\qquad
\m = \m_F + \m'$$
where $\g_Q = (\l_o + \l_F) + \R Z_\D + \m_F$ is the Lie algebra
of the group $G_Q$ and where $\l_o$ is the kernel of non-effectivity 
on $F$, so that $\g_F := Lie(\rho(G_Q)) \simeq \l_F + \R Z_D
+ \m_F$. In \cite{Sp}, the pair $(\g_F, \l_F)$ is called {\it
Morimoto-Nagano pair
of the orbits $G\cdot \eta_t$\/}.\par
\medskip
By the definition of optimal transversal curve, we know that
 there exists a  basis for $\m$, with respect to which 
the complex structure $J_t: \m\to \m$ induced by the CR structure of the orbit
$G\cdot \eta_t$
assumes a particularly simple expression. Using the results in \cite{Sp}, 
one can check that  we may always consider 
a basis for $\R Z_\D + \m$
$$\{ F_0 = Z_\D, F_1, G_1, \dots, F_{n-1}, G_{n-1}\ \}$$
where the following properties hold:  
\roster
\item The elements $F_i$, $G_i$, with $1 \leq i \leq N_F$ are a basis for $\m_F$;
\item the elements $F_i$, $G_i$ with $N_F + 1 \leq i \leq n-1$ are a basis for
$\m'$; 
\item for any $t\neq 0$ (that is when $\eta_t \in \Mreg$),  
$$J_t F_i = 
\left\{\matrix \left[- \tanh^{(-1)^i}(\ell_i t)\right]  G_i & \text{if} \ \
1 \leq i \leq N_F\\
G_i & \text{if} \ \ N_F +1 \leq i \leq n-1
\endmatrix\right.\tag 4.7$$
\endroster
where $\ell_i$ is $2$ if $F_i \in [\m_F, \m_F]^\C \cap \m_F^\C$ and 
 $1$ otherwise. \par
\medskip
There also exists a Cartan subalgebra $\t^\C \subset \l^\C + \C Z_D$ 
for $\g^\C$, so that $\t^\C_F = \t^\C \cap \g^\C_F$ is  a Cartan subalgebra 
for $\g^\C_F$ and the  root systems $R$ and $R_F$ of $(\g^\C, \t^\C)$ and 
$(\g^\C_F, \t^\C_F)$ decompose into
$$R_F = R^o_F + R'_F\ ,\qquad R = R^o \cup R^1 = (R^o_\perp \cup R^o_F) \cup
(R_F'
\cup R')\ ,$$
where 
$$R^o_\perp = \{\ \alpha \ ,\ E_{\alpha} \in \l^\C_o\ \}\quad ,\quad
R^o_F = \{\ \alpha \ ,\ E_{\alpha} \in \l^\C_F\ \}\ ,$$
$$R'_F = \{\ \alpha \ ,\ E_{\alpha} \in \m^\C_F\ \}\quad ,\quad
R' = \{\ \alpha \ ,\ E_{ \alpha} \in \m'{}^\C\}\ .$$
Moreover, the Cartan subalgebra can be chosen so that 
the elements $F_i$, $G_i$ are expressed as follows in 
terms of the root vectors.\par
If $1\leq i \leq N_F$, then two cases may occur: either there exists 
a pair of roots $\{\alpha_i, \alpha^d_i\} \in R'_F$ and an integer $\epsilon_i =
\pm 1$ so that 
$$F_i = \frac{ (E_{\alpha_i} - E_{- \alpha_i}) + 
(-1)^{i+1}\epsilon_i (E_{\alpha^d_i} - E_{- \alpha^d_i})}{2}\ ,\tag 4.8$$
$$G_i = i \frac{(E_{\alpha_i} + E_{- \alpha_i}) + 
(-1)^{i+1}\epsilon_i  (E_{\alpha^d_i} + E_{- \alpha^d_i})}{2}\ ;\tag 4.9$$
or there exists one root $\alpha_i \in R'_F$ such that 
$$F_i = \frac{ E_{\alpha_i} - E_{- \alpha_i}}{\sqrt{2}} \ ,\qquad
 G_i = i \frac{E_{\alpha_i} + E_{- \alpha_i}}{\sqrt{2}}\ .\tag 4.9'$$
This second case may occur only when $N_F$ is odd and, in this case, there
is only one 
such vector $F_i$ and we may assume it is the vector $F_{N_F}$. \par
If  $N_F +1\leq i \leq n-1$, then there exists a root $\beta_i \in R'$ so that 
$$F_i = \frac{ E_{\alpha_i} - E_{- \alpha_i}}{\sqrt{2}} \ ,\qquad
 G_i = i \frac{E_{\alpha_i} + E_{- \alpha_i}}{\sqrt{2}}\ .\tag 4.10$$
The ordering of the roots can be chosen so that the root vectors  $E_{\alpha_i}$ 
which define the elements 
$F_i$, $i\geq N_F +1$, by (4.10) are those corresponding 
to the  roots  $\alpha_i \in R'_+ = R^+ \cap R'$ (see \cite{Sp}).
\bigskip
The explicit expression of the element 
$Z_\D$ for each manifold $M$ of Table 1 is given (up to scaling) 
in the Table of 
Def. 1.7 in \cite{AS}. To determine the  factor to obtain 
the exact expression of 
$Z_\D$, it is possible to consult Table 1 in \cite{Sp},
 where, for each possibility of 
the Morimoto-Nagano pair $(\g_F, \l_F)$, the elements $Z_o$, such that 
$\exp(i t Z_o) \cdot p_o$ is an optimal transversal curve, are given.
For convenience of 
the reader, we list the  expressions 
for $Z_\D$, derived in this way, in Table A1 in Appendix.\par
\bigskip
Let us now write down the explicit expression for the Einstein equations (4.6).
In the following
we will also assume that 
the  Ricci form $\rho$ and the K\"ahler form $\omega$ are {\it tame\/}, i.e. that
the constants 
$C_\rho$ and $C_\omega$ in the expression for 
$Z_\rho(t)$ and $Z_\omega(t)$
are both equal to $0$, even in case $M$ is a non-standard KE-manifold. \par
We will see a posteriori that, for all the
K-manifold of considered in Main Theorem (2), the Ricci form 
is tame; hence  we have no loss of generality  with this assumption.\par
\medskip
By looking at Table 1 in \cite{Sp},  one can check that for any $1 \leq i \leq
N_F$, the 
bracket $[F_{i}, G_{i}]$ is orthogonal to $\z(\l)$. Set $i = 1$. Then, 
since $I_\omega, I_\rho \in \z(\l)$, it follows that 
$(4.6)_1$ holds if and only if  for any $t\in \R$
$$\B(Z_\rho(t), [F_1, G_1]) - c \B(Z_\omega(t), [F_1, G_1]) = 
\frac{\rho_{\eta_t}(\hat F_1, J \hat F_1) -
 c \omega_{\eta_t}(\hat F_1, J \hat F_1)}{- \coth(\ell_1 t)} = 0\ .\tag 4.11$$
Using the expression for $\rho_{\eta_t}$ 
given in Thm. 5.1 and  Prop. 5.2 in \cite{Sp}, assuming that 
$F_1$ is chosen so that 
$[F_1, G_1]_{\l+ \m}$ is orthogonal to $[\l + \m, \l + \m]$ 
(there is no loss of generality 
in this assumption) and by (5.13) and (5.13') in \cite{Sp},
it follows that  (4.11) is equal to 
(we used the fact that 
for any $1 \leq j \leq N_F$, $\B(Z_\D, [F_j, G_j]) = \ell_j$ - 
see (5.13) in \cite{Sp})
$$ h'(t)  -
2  \sum_{i = 1}^{N_F} \ell_i \tanh^{(-1)^{i+1}}(\ell_i t)
   +
\tilde c  f_\omega(t)  + 4 \B(Z^{\kappa}, Z_\D) = 0\
,\tag 4.12$$
where 
$$h(t) = \log \det [g_{ij}(t)]\ ,\qquad 
g_{ij}(t) = g_{\eta_t}(\hat F_i , \hat F_j) = - \omega_{\eta_t}(\hat F_i, J \hat
F_j)\ ,\tag 4.13$$
$$\tilde c =  - 2\B(Z_\D, Z_\D) c >0\ ,\qquad
Z^{\kappa} = \sum_{k = N_F +1}^{n-1} i H_{\beta_k}\ .\tag4.14$$
\par
From the properties of the chosen
adapted basis, it follow that the vectors $\hat F_i|_{\eta_t}$, $J \hat
F_i|_{\eta_t}$ are $g$-orthogonal 
for any $t \neq 0$ and that  for  $1\leq 2s, 2s+1\leq N_F$ and  $N_F +1 \leq
i\leq n-1$:
$$g_{\eta_t}(\hat F_0, \hat F_0) =  - \B(Z_\D, Z_\D) f'_\omega(t)\ ,
\tag4.15$$
$$
g_{\eta_t}(\hat F_{2s +1}, \hat F_{2s +1}) = \coth(\ell_{2s+1} t) f_\omega(t)
\ell_{2s+1}\ ,\  g_{\eta_t}(\hat F_{2s}, \hat F_{2s}) = \tanh(\ell_{2s} t)
f_\omega(t) \ell_{2s}\ ,\tag4.16$$
$$g_{\eta_t}(\hat F_i, \hat F_i) = - \beta_i(i Z_\D)f_\omega(t) - \beta_i(i
I_\omega) \ .
\tag4.17$$
In case $\beta_i(Z_\D) \neq 0$ (which is always true if the regular 
$G$-orbits are Levi non-degenerate), we may write 
$$g_{\eta_t}(\hat F_i, \hat F_i) = - \beta_i(i Z_\D)(f_\omega(t) + a_{i})\
,\qquad
\text{where}\quad a_{i} = \frac{\B(i H_{\beta_i}, I_\omega)}{
\B(i H_{\beta_i}, Z_\D)} \ .
\tag4.17'$$
Now,  if $N_F$ is odd,  
$$h'(t) = \frac{f''_\omega(t)}{f'_\omega(t)} + 
N_F\frac{f'_\omega(t)}{f_\omega(t)}
+  
f_\omega'(t) \sum_{i = N_F+1}^{n-1} \frac{1}
{f_\omega(t) + a_i} +
\ell_{N_F}(\tanh(\ell_{N_F} t) - \coth(\ell_{N_F} t))\ .$$
In case  $N_F$ is even, we have 
$$h'(t) = \frac{f''_\omega(t)}{f'_\omega(t)} + 
N_F\frac{f'_\omega(t)}{f_\omega(t)}
+  
f_\omega'(t) \sum_{i = N_F+1}^{n-1} \frac{1}
{f_\omega(t) + a_i}\ .$$
So,  (4.12)  becomes
$$
\boxed{
\matrix
\ \\
\frac{f''_\omega(t)}{f'_\omega(t)} + f'_\omega(t) \left(\frac{ N_F
}{f_\omega(t)} + 
\sum_{j = N_F+1}^{n-1} \frac{1}
{f_\omega(t) + a_j}\right) + \tilde c f_\omega(t) - \\
\ \\
- N^{(1)}_F (\tanh(t) +\coth(t)) - 2 N^{(2)}_F (\tanh(2t) +\coth(2t))  + 
4 \B(Z^{\kappa}, Z_\D) = 0\\
\ 
\endmatrix
} 
\tag 4.18$$
where $N^{(1)}_F$ denotes the number of vectors $F_i$ with $\ell_i = 1$ and 
$N^{(2)}_F$ the number of those vectors with $\ell_i = 2$.\par
\medskip
Let us now consider the condition $(4.6)_2$.
Recall that  $\g_Q = (\l_o + \l_F) + \R Z_\D + \m_F$ 
is the isotropy subalgebra of the flag manifold $G/G_Q$ and 
that the restrictions $\beta_m|_{\z(\l)}$ of the roots
$\beta_m \in R'$, corresponding to  highest weight vectors of $\m'{}^\C$,
generate
the entire dual space $\z(\l)^\C{}^*$
(see e.g. \cite{Al}). \par
There is no loss of generality if we assume that for any such roots $\beta_m$,
the vector
 $F_m = \frac{E_{\beta_m} - E_{- \beta_m}}{\sqrt{2}}$ coincides with one 
of the vectors $F_j$, $N_F +1 \leq j \leq n-1$, which belong to the basis for $\m
= \m_F + \m'$
we have chosen.\par
Consider one of these vectors $F_m$ and the equation 
$$\rho_{\eta_{t_o}}(\hat F_m, J \hat F_m) - c \omega_{\eta_{t_o}}(\hat F_m, J\hat
F_m) = 0\ .\tag 4.19$$
Note that equation (4.19) is equivalent to any other 
equation, which is obtained by  replacing $F_m$ by 
any other vector $F_i$ in the same irreducible $\l$-module in $\m'$ (which is
clearly included in some
irreducible $\g_Q$-module).\par 
\medskip
From Thm. 5.1 and Prop. 5.2 in \cite{Sp}, (4.19) is equivalent to
$$
   h'(t)  - 
 2 \sum_{i = 1}^{N_F} \ell_i \tanh^{(-1)^{i+1}}(\ell_i t) 
+ \tilde c f  +  2 \B(Z^{\kappa}, Z_\D)
 + $$
$$ +
\frac{2 \B(Z_\D, Z_\D) \B(Z^{\kappa}, i H_{\beta_m})}{\B(Z_\D, i H_{\beta_m})}
- \frac{ 2c \B(Z_\D, Z_\D) \B( I_\omega, i H_{\beta_m})}
{\B(Z_\D, i H_{\beta_m})}   = 0\ , \tag 4.20$$
Subtracting (4.12) from (4.20), we get 
that 
$$0 =  -  2 \B(Z^{\kappa}, Z_\D)+
\frac{2 \B(Z_\D, Z_\D) \B(Z^{\kappa}, i H_{\beta_m})}{\B(Z_\D, i H_{\beta_m})}
- \frac{ 2c \B(Z_\D, Z_\D) \B( I_\omega, i H_{\beta_m})}
{\B(Z_\D, i H_{\beta_m})}\ ,$$
so that
$$\B( c I_\omega, i H_{\beta_m}) = 
\B(Z^{\kappa} - \frac{\B(Z^{\kappa}, Z_\D)}
{\B(Z_\D, Z_\D)} Z_\D ,i H_{\beta_m}) \tag 4.21$$
From the previous remarks,  
it follows that (4.21) 
determines 
$I_\omega$ uniquely and we may write
$$\boxed{
I_\omega =  \frac{1}{c} \left( Z^{\kappa} - \frac{\B(Z^{\kappa}, Z_\D)}
{\B(Z_\D, Z_\D)} Z_\D \right) = 
\frac{1}{c} Z^{\kappa}_{\perp}
}\tag 4.22$$
where by $Z^{\kappa}_{\perp}$ we denote the orthogonal 
projection of $Z^{\kappa}$ in $(Z_\D)^\perp$. \par
It also implies that the coefficients $a_{m}$, which appear in (4.17'), are equal
to
$$\boxed{
a_{m} =  \frac{1}{c} \frac{\B(i H_{\beta_m}, Z^{\kappa}_{\perp})}{\B(i
H_{\beta_m}, Z_\D)}
}\tag
4.23$$
Equation (4.18) and (4.22) are the Einstein equations we were looking for.\par
\medskip
The explicit  expressions for $Z^\kappa$ and the values of the coefficient
$a_{m}$
for all K-manifolds of Table 1 are listed in Table A1 (see Appendix).  
As remarked in Appendix, 
for all those cases $Z^\kappa = Z^\kappa_\perp$ and hence {\it 
(4.22) reduces 
to $I_\omega = \frac{1}{c} Z^\kappa$\/}.
\bigskip
\bigskip
\subsubhead 4.3 The differential problem which characterizes the
Einstein-K\"ahler metrics
\endsubsubhead
\bigskip
We are now ready for the main result of this section. The following Theorem gives
the differential 
problem that one has to solve in order to determine the K\"ahler-Einstein metrics
(if any) 
on the  non-standard K-manifolds of Table 1. \par
\bigskip
In all the statements of this subsection, we will assume the same hypothesis 
on $\tilde M$ and the same optimal transversal curve $\eta_t$ used in 
\S 4.2.\par
\medskip
\proclaim{Theorem 4.2} Let $\tilde M = G\times_{G_Q, \rho} F$ be one 
of the manifolds given in Table 1.  Then $\tilde M$ 
admits a K\"ahler-Einstein metric with Einstein
constant $c>0$ and tame Ricci form $\rho$ if there exists a smooth 
function $f: ]0, +\infty[ \subset \R \to ]0,+ \infty[$ 
which verifies the following conditions:
\roster
\item for any $t\in ]0, \infty[$ and  for  any root $\beta_m \in R'_+$
$$f(t) >0\ ,\quad f'(t) >0\ ,\quad - (f(t) + a_{m}) \beta_m(i Z_\D)>0\ ,$$
where the coefficient  $a_{m}$  are defined by (4.23);
\item $f$ verifies on $(0, \infty)$ the  differential equation (4.18) with 
$\tilde c$ and $Z^\kappa$ as defined in (4.14);
\item $- \B(i H_{\beta_m}, I_\omega) = 
 - a_{m} \beta_m(i Z_\D) > 0$ for any $\beta_m \in R'_+$, 
$\lim_{t\to 0} f(t) = 0 =\lim_{t\to 0}f''(t)$ and the following limits 
exist and are finite
$$\lim_{t\to 0} f'(t) = C_1\ ,\qquad \lim_{t\to 0}f'''(t) = C_2$$
with $C_1 > 0$; 
\item the limits 
$\lim_{t\to + \infty} f(t)$ and $\lim_{t\to + \infty} - (f(t) + a_m) \beta_m(Z_\D)$
are finite and positive for any $\beta_m \in R'_+$ and 
$$\lim_{t \to + \infty} e^{2 \varepsilon_F t} f'(t) = C_3\ ,\quad
\lim_{t \to + \infty} e^{\varepsilon_F t} \left(1 + \frac{1}{\varepsilon_F}
\frac{f''(t)}{2
f'(t)} \right) = 0\ , $$
$$ \lim_{t \to + \infty} e^{2 \varepsilon_F t} \left(1 + 
\frac{1}{\varepsilon_F} \frac{5f''(t)}{6
f'(t)}  + \frac{1}{\varepsilon_F^2}\frac{f'''(t)}{6
f'(t)}\right) = C_4\ ,$$
for some finite values $C_3 >0 $, $C_4$ and 
where   $\varepsilon_F = \left\{\matrix 2 & \text{if}\ \  F = \CP^r\hfill\\
1 & \text{if}\ \  F = Q^r\hfill\endmatrix\right.$.
\endroster
If (1) to (4) are satisfied, then the K\"ahler form $\omega$ 
of the K\"ahler-Einstein metric is given  by the algebraic
representative 
$$Z_\omega(t) = f(t) Z + {1\over c}Z^\kappa\ ,$$
where $f$ has to be meant as extended  over the whole real axis as 
a continuous odd function. 
\endproclaim
\demo{Proof}  We claim  that   a $G$-invariant metric $g$ on $\Mreg (= \tilde
\Mreg)$
with tame K\"ahler form $\omega$ 
 is  K\"ahler-Einstein if and only if
the algebraic representative $Z_\omega(t) = f(t) Z_\D + I_\omega$, with $t > 0$,
verifies (1), (2) and $I_\omega = {1\over c}Z^\kappa$. In fact, by the 
previous discussion,  
we know that (4.18) and $I_\omega = {1\over c}Z^\kappa$ are 
necessary and 
sufficient conditions for the algebraic representative
$Z_\omega(t) = f(t) Z_\D + I_\omega$ to represent a K\"ahler-Einstein metric on 
$\Mreg$. 
Furthermore, by (4.15), (4.16) and (4.17),  $g$ is  
positive definite on $\Mreg$ if and only if the conditions (1) are verified.\par
To check that $f$ has to be extended as an odd function over $\R$, we recall 
that, in case $\tilde M$ is a KO-manifold, 
$\eta$ coincides with a re-parameterization 
of a  normal geodesic for some $G$-invariant K\"ahler metric (see \cite{Sp}).
Let  $\gamma_s$ be one of these normal geodesics,  with 
$$\gamma_0 = \eta_0 = p_o\ ,\qquad \gamma_s = \eta_{t(s)}\ ,
\qquad \frac{ds}{dt} =
||\hat Z_\D||_{\eta_t}.$$
By definition of geodesic symmetry $\sigma$ at $p_o$ (see
 \cite{AA}, \cite{AA1} for the definition)
we have that $\sigma_*(\hat Z_\D|_{p_o}) = 
-\hat Z_\D|_{p_o}$. Then
$$f(t(-s))  Z_\D + I_\omega = 
Z_{\omega}(t(-s)) = Ad(\sigma_{p_o})Z_{\omega}(t(s)) = -  f(t(s)) Z_\D +
I_\omega\ ,\tag 4.25$$
and therefore $f(t(-s)) = - f(t(s))$. Now, by considering that
 ${dt\over ds} = \frac{1}{||\hat Z_\D||_{\eta(t(s))}}$
is an even function of $s$  (and hence that  $t(s)$ is an odd function), it
follows that 
$f:\R \to \R$  is an  odd function. When $\tilde M$ is a KE-manifold, 
we observe that $\eta$ lies in $Fix(L) = \CP^1\times \CP^1$ or $\CP^2$
and we may argue in a similar way, replacing $\tilde M$ with 
$Fix(L)$. \par
\medskip
 By a result of De Turk and Kazdan
 (see e.g. \cite{Be}), any $\Cal C^2$ Einstein metric 
 is   real analytic in a geodesic normal coordinate 
system. Therefore it remains to prove that, if $f: ]0, + \infty[ \to \R$, 
  verifies 
(1) and (2),
then the K\"ahler-Einstein metric  $g$ on
$\Mreg$, associated 
with  $Z_\omega = f(t) Z_\D + {1\over c}Z^\kappa$, for $t>0$, 
 extends to a $\Cal C^2$ metric on $\tilde M$ 
 if (3) and (4) hold.
\smallskip
A key ingredient to prove the $\Cal C^2$-extendibility at the noncomplex
singular orbit
is given by the following Lemma.\par
\enddemo
\medskip
\proclaim{Lemma 4.3} Let  $g$ be the K\"ahler-Einstein metric on $\Mreg$
associated with the K\"ahler form $\omega$
given by the algebraic representative $Z_\omega(t) = f(t) Z_\D + I_\omega$, 
with $I_\omega = Z^\kappa$ and $f : \R \to \R $ an odd function which is 
smooth on $\R \setminus\{0\}$ and  verifies
(1) and (2) of Theorem 4.2. Define
$\Lambda: ]0, + \infty[ \to \R$ as
$$\Lambda(t) = - \B(Z_\D, Z_\D)\int_0^t f(u) du\ .\tag 4.26$$
Let $\Cal U$ be a $G$-invariant neighborhood of 
the noncomplex singular $G$-orbit in $\tilde M$. 
If the smooth tensor field $dd^c \hat \Lambda$ where
$$\hat \Lambda: \Mreg \cap \Cal U\to \R\ ,\qquad \hat \Lambda(p) = 
\hat \Lambda(g\cdot \eta_t) = \Lambda(t)\ ,\tag 4.27$$
extends as a $\Cal C^2$ tensor field on the whole $\Cal U$, then
$\omega$ extends to a
$\Cal C^2$ 2-form on $\Cal U$. In particular $\omega$ extends if 
$\hat \Lambda$ extends 
as a $C^4$ function.
\endproclaim
\demo{Proof} Consider the unique 
$G$-invariant $J$-invariant closed 2-forms $\varpi_f$ and 
$\varpi_{I_\omega}$ on $\Mreg$ 
with associated algebraic representatives $Z_{\varpi_f}(t) \= f(t) Z_\D$ and 
$Z_{\varpi_{I_\omega}} \= I_\omega$. From definitions, 
the K\"ahler form $\omega$ on $\Mreg$,  determined by $Z_\omega$, coincides with 
$$\omega = \varpi_f + \varpi_{I_\omega}\ .$$
We claim that $\varpi_{I_\omega}$ extends to a 
smooth $G$-invariant 2-form on $\tilde M$. In fact,
we recall that $\frac{1}{c} Z^\kappa$ is the element in the center 
$Z(G_Q)$, which corresponds to
the $G$-invariant K\"ahler-Einstein metric $\hat g$ on the flag manifold $G/G_Q$,
with 
K\"ahler form defined by (see \cite{Be}, Ch.8)
$$\hat \omega_{o}(\hat X', \hat Y') = \frac{1}{c}\B(Z^\kappa, [X', Y'])\ ,$$
where $o = eG$ and for any $X', Y'\in \m'$. So, by  definitions and Proposition
4.1, 
at any point $\eta_t \in \Mreg$ and for any  vectors 
$X, Y\in  \R Z_\D + \m_F \subset \g_Q$ and  
$X', Y' \in \m'$
$$\varpi_{I_\omega}(\hat X, \hat Y)_{\eta_t} = \frac{1}{c}\B(Z^\kappa, [X,Y]) =
0\ ,
\varpi_{I_\omega}(\hat X, \hat X')_{\eta_t} = \frac{1}{c}\B(Z^\kappa, [X,X']) =
0\ ,$$
$$\varpi_{I_\omega}(\hat Z_\D, J \hat Z_\D)_{\eta_t} = 0\ , 
\varpi_{I_\omega}(\hat X', J \hat Y')_{\eta_t} = 
\pi^*(\hat \omega)(\hat X', \hat Y')_{\eta_t}\ ,$$
where $\pi$ is the  projection $\pi: \tilde M \to G/G_Q$. This implies that 
on $\Mreg$,
$\varpi_{I_\omega} = \pi^*(\hat \omega)$ and this proves that $\varpi_{I_\omega}$
can be extended smoothly on the entire $\tilde M$. Hence,
 $\omega$ extends to a $\Cal C^2$ 2-form on $\Cal U$ if 
and only if $\varpi_f$ does. \par
Now, notice that at any point $p \in \Mreg$ and for any 
two vector fields $V, W$
$$\varpi_f(V,W)|_p = 
(d d^c \hat \Lambda)(V,W)|_p \= \left[- V(JW(\hat \Lambda)) + 
W(JV(\hat \Lambda)) + J[V, W](\hat \Lambda)\right]_p\ .\tag 4.28$$ 
Due to $G$-invariance, (4.28) needs to be checked only 
at the points 
$\eta_t \in \Mreg$, with vector fields $V, W$ of the form 
$\hat Z_D$, $J \hat Z_\D$ or 
$\hat X, \hat Y$ with $X, Y \in \m$. Since $\hat \Lambda$ 
is constant along  $G$-orbits, it follows that 
$$(d d^c \hat \Lambda)(\hat X,\hat Y)|_{\eta_t} = 
- J\widehat{[X, Y]}(\hat \Lambda)_{\eta_t} = 
-  \frac{\B(Z_\D, [X,Y])}{\B(Z_\D, Z_\D)} J \hat Z_\D(\hat \Lambda)_{\eta_t} = $$
$$ = \B(f(t) Z_\D, [X,Y]) = 
\varpi_f(\hat X,\hat Y)|_{\eta_t}$$
$$(d d^c \hat \Lambda)(\hat X,\hat Z_\D)|_{\eta_t} = 
- J\widehat{[X, Z_\D]}(\hat \Lambda)_{\eta_t} = 0 = 
\varpi_f(\hat X,\hat Z_\D)|_{\eta_t}\ ,$$
$$
(d d^c \hat \Lambda)(\hat X,J \hat Z_\D)|_{\eta_t} = 
J\widehat{[J_t X, Z_\D]}(\hat \Lambda)_{\eta_t} = 0 = 
\varpi_f(\hat X,J \hat Z_\D)|_{\eta_t}\ ,$$
$$(d d^c \hat \Lambda)(\hat Z_\D, J \hat Z_\D)|_{\eta_t} =  
J \hat Z_\D(J \hat Z_\D(\hat \Lambda)|_{\eta_t} = 
- \B(Z_\D, Z_\D) f'(t) = \varpi_f(\hat Z_\D, J\hat Z_\D)|_{\eta_t}\ ,$$
and this proves (4.28). Our claim follows immediately.\qed
\enddemo
\bigskip
\demo{End of proof of Theorem 4.2} Let us now prove that (3) is a sufficient
condition for the  K\"ahler-Einstein metric $g$  on $\Mreg$ to be 
$\Cal C^2$ extendible on a neighborhood of the noncomplex singular orbit
$G\cdot \eta_0$.\par 
First of all, we claim that we can find a slice $S\ni \eta_0$, containing 
$\eta_t$ for small $t$. Indeed, if $\tilde M$ is 
KO-manifold, then $\eta_t$ is reparameterization of
a normal geodesic for any $G$-invariant K\"ahler metric
and therefore one can take $S$ to be $\exp_{\eta_0}(V)$
for some suitable open neighborhood of $0$ in the normal 
space to the singular orbit $G\cdot \eta_0$; when $\tilde M$ is 
a KE-manifold, then $S$ can be constructed along the same line 
inside the fixed point set $F \subset Fix(L)$.\par
We also know that 
there exists a local section $\chi: \Cal V_{\eta_0} \subset G/G_{\eta_0} \to G$
such that the map $(v,s) \to \chi(v)\cdot s$ is a 
 diffeomorphism between $\Cal V \times S$ and an open neighborhood of 
$\eta_0$ in $G\cdot S \subset \tilde M$ (see Lemma 2.2 in \cite{Pa}).\par
From this and from the fact that the function 
$\hat \Lambda$, defined  in (4.27), is constant along  $G$-orbits, 
it follows that it extends as a $\Cal C^4$ function on
a neighborhood of $\eta_0$
if and 
only if 
$\hat \Lambda|_{S\setminus \{\eta_0\}}: S\setminus \{\eta_0\} 
\to \R$ extends as a $\Cal C^4$ 
 function over $S$.  
\par
We identify $S$ with a suitable ball of radius $r$ in $\R^{2n - N_F - 1}$, 
via some $G_{\eta_0}$-equivariant identification. Let us also denote by 
$T: ]-r, r[ \to \R$
 the odd smooth function  defined on the positive values by the relation 
$$g\cdot x = \eta_{T(|x|)}\ ,\quad \text{for}Ý\ x\in S\ ,\ g\in G_{\eta_0}\ .$$
Then  
$$\hat \Lambda: S \setminus \{0\} \to \R\ ,\quad \hat \Lambda(x) =
\Lambda(T(|x|))$$
 by the $G_{\eta_0}$-invariance. So, 
necessary and sufficient condition for $\hat \Lambda$ to extend to a $\Cal C^4$
function on $S$ is that $\Lambda: ]0, r[\to  ]0, + \infty[$ extends to 
a $\Cal C^4$ even function over $]-r,r[$. 
This occurs if 
(3) (without the restriction on the sign of $C_1$) is verified. 
Moreover, using (4.15) - (4.17), 
one can check that if  (1), (2) and (3) holds with $C_1>0$, 
then the extension of 
$g$ on  
$G\cdot \eta_0$  is positive definite 
 at $\eta_0$ and hence on the whole orbit.\par 
\medskip
We now deal with the extendability of the metric on a 
$G$-invariant neighborhood of the 
complex singular orbit $G\cdot p_o$, where $p_o = \lim_{t\to \infty} \eta_t$.
Note 
that $G_{p_o}$ is connected and its Lie algebra is given by 
$\g_{p_o} = \l + \R\cdot Z_\D$. \par
As before, we consider a 2-dimensional slice $S$ through 
$p_o$, which contains the curve 
$\eta_t$ for large $t$, and, using the Lemma 2.2 in \cite{Pa}, 
we determine a  
neighborhood $U$ of $p_o$ in $G\cdot p_o$ and a local section $\chi:U \to G$ 
such that the map $U\times S \ni (x,s)\mapsto \chi(x)\cdot s$ is a diffeomorphism
onto its image.\par 
Now notice that, since  $S$ is $G_{p_o}$-invariant, 
the vector field $\hat Z_{\D}$ is tangent to $S$. So 
any
tangent space $T_{\eta_t} S$ is spanned by the vectors 
$\hat Z_D|_{\eta_t}$ and 
$J\hat Z_\D|_{\eta_t} = \eta'_t$ and hence it is $J$-invariant; in particular, 
$S$ is a $G_{p_o}$-invariant complex submanifold of 
$\tilde M$. \par
By the Riemann Mapping Theorem, by choosing $S$ sufficiently small, 
we may assume that there exists a biholomorphism 
$\varphi: S \to \Delta \subset \C$ which maps $S$ onto the unit disc 
$\Delta = \{ |z| < 1\} \subset \C$ and so that $\varphi(p_o) = 0$. 
If we use $\varphi$ to identify $S$ with $\Delta$, it follows from 
Schwarz Lemma  
that the 1-parameter group $\exp(\R Z_\D) \subset G_{p_o}$ acts on 
$\Delta$ as a closed group of rotations.  In particular,
using the standard 
polar coordinates $(r,\theta)$  for $\Delta$,  we have that
$$ \hat  Z_\D|_{S} = k\frac{\partial}{\partial \theta}\ ,\quad
J\hat Z_{\D}|_S = - k r \frac{\partial}{\partial r}\ ,\quad
\eta_t = (r(t) = A e^{ - k t}, \theta = 0)\ ,\tag 4.29$$
for some positive real constant $A$ and 
where $\frac{1}{k}$ is the smallest positive rational number such that 
$\exp(\frac{2\pi}{k} Z_\D) \cdot \eta_t = \eta_t$ for any $t>0$. \par
Let us now prove that $k = \varepsilon_F$ as defined in the 
statement of the Theorem. Using Table 1 in \cite{Sp}, 
where the explicit expressions of $\rho_*(Z_\D)$ 
are given, one can check that the adjoint action of 
$\rho(\exp(\pi Z_\D))$ on the tangent space 
$T_{p_o} (\rho(G_Q)\cdot p_o)$
of the singular orbit $G\cdot p_o \subset F$ is equal to  $- I$
and that for no value $0 < k <1$, $\rho(\exp(\pi k Z_\D))$
acts on $T_{p_o} (\rho(G_Q)\cdot p_o)$ as a multiple of the identity. 
On the other hand, by direct inspection of the 
explicit action of   
$\rho(G_Q) = \SO_{r+1}, \SO_{r+1}/\Z_2$ on $F = Q^r, \CP^r$
(see e.g. \cite{Uc}, p. 157-158), 
the action of $\rho(\exp(\pi Z_\D))$ on 
$S$ coincides  with the geodesic symmetry 
if $F = Q^r$ and  with  
identity map 
in case $F = \CP^r$. In particular, it follows that 
 $k = \varepsilon_F$.\par
\medskip
From the fact that $||\hat Z_\D||^2_{\eta_t} = 
- \B(Z_\D, Z_\D) f'(t)$,  we get that the
restriction of $g$ on the vectors tangent to  $S \setminus \{p_o\}$ 
is  
$$g|_{TS
 \times  TS }= \frac{- B f'(t(r))}{\varepsilon^2_F r^2} 
(dr^2 + r^2 d\theta^2)\ .\tag 4.30$$
By the results in \cite{Ve} (see also \cite{Ber}), we have that
$g|_{TS
 \times  TS }$ extends to a $\Cal C^2$ Riemannian
metric  over the whole $S$ if and only if the function 
$\frac{ \sqrt{f'(t(r))}}{r}$ extends as a $\Cal C^2$ even positive 
function 
on  $]-1, 1[$. This is equivalent to  conditions (4).\par
It remains to prove that if (4) is verified then the metric $g$ extends 
to  $\Cal C^2$ Riemannian metric over the entire $\chi(U) \cdot S \simeq 
U \times S$.
\medskip
Consider now the following frames field 
 $\{X_1, \dots, X_n\}$.  Let  
$X_1, X_2$ be any two smooth vector fields
tangent to $S$ and linearly independent
at any point. Then, consider  the vectors  $E_i\in \g$, $3 \leq i \leq 
2n$, which are equal to the linearly independent vectors $F_j$ or $G_j$
defined in (4.8) - (4.10),  and for any $x\in S$
 let $X_{i}|_x = \hat E_i|_x$. By choosing 
$S$ sufficiently small, we may always suppose that 
all the $X_i$'s are linearly independent at any point of 
$S$.  Finally,  let us extend the field of frames $\{X_i\}$
on  $S$ 
to a field of frames on $\chi(U) \cdot S \simeq U \times S$, 
by setting 
$$\{X_1, \dots, X_n\}|_{\chi(u)\cdot x} =  \{\chi(u)_*(X_1|_x),
 \dots, \chi(u)_*(X_1|_x)\}\ .$$
Now,  $g$ extends to a $\Cal C^2$ tensor if and only if
the functions $g_{ij} = g(X_i, X_j)$ extends as $\Cal C^2$ functions. 
By the $G$-invariance of $g$ and by the construction of the 
frames $\{X_i\}$,  these conditions are equivalent to the conditions 
for the extendability  of the restrictions
$g_{ij}|_{S\setminus \{p_o\}}$ to $\Cal C^2$ functions  over $S$. \par
Observe that, 
by definitions from the fact that  $[Z_\D, Z_\omega(t)] = 0$ for any $t$, 
we have for $i,j\geq 3$ and for any point 
$x = (r, \theta) \in S\setminus \{p_o\}$
$$g_{ij}|_{(r, \theta)} = \left.
\left[\exp(-\frac{\theta}{\varepsilon_F} Z_\D)^*g\right]
(\hat E_{i}, \hat E_{j})\right|_{(r,\theta)} =  $$
$$ =  - 
\B\left(Z_\omega(t(r)), [Ad_{\exp(-\frac{\theta}{\varepsilon_F} Z_\D)}(E_i), 
J_t Ad_{\exp(-\frac{\theta}{\varepsilon_F}  Z_\D)} (E_j)]\right) =$$
$$ = 
- \B\left(Ad_{\exp(\frac{\theta}{\varepsilon_F}  Z_\D)}
\left(Z_\omega(t(r))\right), [E_i, 
 J^\theta_t E_j]\right) =  - \B\left(Z_\omega(t(r)), [E_i, 
 J^\theta_t E_j]\right)\ ,\tag 4.31$$
where
$$J^\theta_t = Ad_{\exp(\frac{\theta}{\varepsilon_F}  Z_\D)}\circ J_t 
\circ Ad_{\exp(-\frac{\theta}{\varepsilon_F}Z_\D)}\ .\tag 4.32$$
From the the fact that the projection  $\pi: \tilde M \to G/G_Q$ is holomorphic, 
it follows that $\operatorname{ad}_{Z_\D} \circ J_t|_{\m'} = 
J_t \circ \operatorname{ad}_{Z_\D}|_{\m'}$ for any $t$ and hence 
that $J^\theta_t|_{\m'} \equiv J_t$; using the explicit 
expressions of  $Z_\D$ given in Table 1 in \cite{Sp},
one can check that, for any $1 \leq i\leq N_F$ 
$$Ad_{\exp(\frac{\theta}{\varepsilon_F}Z_\D)}(F_i)
 = \cos\left(
\frac{ \theta}{\varepsilon_F}\right) F_i -
\sin\left(
\frac{ \theta}{\varepsilon_F}\right) G_i\ ,$$
$$
Ad_{\exp(\frac{\theta}{\varepsilon_F}Z_\D)}(G_i)
 =  \sin\left(
\frac{ \theta}{\varepsilon_F}\right)F_i +
\cos\left(
\frac{ \theta}{\varepsilon_F}\right) G_i\ ,$$
and hence that 
$$J^\theta_t(F_i) = 
(-1)^{i}\frac{\tanh(t) - \coth(t)}{2} \sin\left(
\frac{2 \theta}{\varepsilon_F}\right) F_i - $$
$$ - 
\left(\tanh^{(-1)^i}(t) \cos^2
\left(
\frac{\theta}{\varepsilon_F}\right) + 
\tanh^{(-1)^{i+1}}(t) \sin^2
\left(
\frac{\theta}{\varepsilon_F}\right)
\right) G_i \ ,$$
$$J^\theta_t(G_i) = 
\left(\tanh^{(-1)^i}(t) \cos^2
\left(
\frac{\theta}{\varepsilon_F}\right) + 
\tanh^{(-1)^{i+1}}(t) \sin^2
\left(
\frac{\theta}{\varepsilon_F}\right)
\right) F_i
-  $$
$$ - (-1)^{i}\frac{\tanh(t) - \coth(t)}{2} \sin\left(
\frac{2 \theta}{\varepsilon_F}\right) G_i \ .$$
From this, (4.15) and (4.16),
we conclude that the functions $g_{ij}(r, \theta)$, with $3 \leq i,j \leq 2 N_F$
and  which are not 
identically equal to $0$,  are proportional to functions of the form
$$f(t(r)) + a_i\ ,\qquad \frac{\tanh(t(r)) - \coth(t(r))}{2} \sin\left(
\frac{2 \theta}{\varepsilon_F}\right) f(t(r))\ ,\tag 4.33$$
$$
\left(\tanh^{(-1)^i}(t(r)) \cos^2
\left(
\frac{\theta}{\varepsilon_F}\right) +
\tanh^{(-1)^{i+1}}(t(r)) \sin^2
\left(
\frac{\theta}{\varepsilon_F}\right)
\right) f(t(r))\ .\tag 4.34$$
Now, by the proof of 
Prop. 2.7 in \cite{KW}, the functions $g_{ij}$ of the form 
(4.33) and (4.34) extend as $\Cal C^2$ functions on the origin if and
only if 
\roster 
\item"i)" the value $\lim_{r \to 0} g_{ij}(r, \theta)$ does not
depend on $\theta$;
\item"ii)" $g_{ij}$ extends as a $\Cal C^2$ function 
over $]-1, 1[ \times \R$ such that
$g_{ij}(- r, \theta) = g_{ij}(r, \theta + \pi)$ for any 
$\theta$ and any $r>0$; 
\item"iii)" for $k = 1,2$, 
 $\lim_{r \to 0} r^k \frac{\partial^k g_{ij}(r, \theta)}
{\partial r^k}$ is a homogeneous 
polynomial of degree $k$ in the variables
$x = r \cos(\theta)$ and $y = r \sin(\theta)$.
\endroster
First of all, note that 
(i) is verified by all functions (4.33), (4.34).
Moreover, by some straightforward computations, 
one can check
that all those functions verify (ii) if and only if 
$f(t(r))$ extends as a $\Cal C^2$ even function on a neighborhood of 
$0$, which turns out to be equivalent to the conditions 
$$\lim_{t\to \infty} f(t) = C_5 \ ,\quad 
 \lim_{t\to \infty} e^{\varepsilon_F t} f'(t) = 0\ ,\quad 
\lim_{t\to \infty} e^{2 \varepsilon _F t} \left(f'(t) + 
\frac{f''(t)}{\varepsilon_F}\right) = C_6$$
for some finite values $C_5$ and $C_6$; 
but these conditions are immediately satisfied if $f$ verifies (4).
In a similar way, a simple computation shows that also (iii) is verified 
by all functions (4.33), (4.34) since $f$ verifies (4). 
\par
Finally,  we recall that, by the previous remarks,
if conditions (4) are satisfied, 
then the functions
$g_{ij}|_{S\setminus \{p_o\}}$ with $i,j = 1,2$ are immediately extendible 
on $S$, while
the functions $g_{i\ell}|_{S\setminus \{p_o\}}$ 
with $i = 1,2$ and $\ell \geq 3$ are
identically vanishing. \par
The positivity of the extension of the metric on the points of the 
singular complex orbit is assured by the positivity of the limits 
$\lim_{t\to + \infty} f(t)$ and 
$\lim_{t\to + \infty} - (f(t) + a_m) \beta_m(Z_\D)$.
 \qed
\enddemo
\bigskip
\bigskip
\subhead 5. New examples of non-standard K\"ahler-Einstein K-manifolds
\endsubhead
\bigskip
We have now all ingredients to determine new examples of
 non-standard K-manifolds that belong to the
class described by Corollary 3.5 (ii). The  
main result is the following.\par
\bigskip
\proclaim{Theorem 5.1} Let $\tilde M = G\times_{G_Q, \rho} F$ be one 
of the manifolds given in Table 1, but 
 not in case 1 with $G = \SU_3$
and $F = \CP^2$, not in case 2 with $G = \SU_p \times \SU_2$, $p>2$, not
in case 4 with $F = Q^7$, nor in case 5 with $F = Q^9$. 
Then it is K\"ahler-Einstein with positive first Chern class.
\endproclaim
\medskip
The crucial step for the proof of Theorem 5.1 is the following proposition.
The meaning of any adopted notation is the same of \S 4.\par
\medskip
\proclaim{Proposition 5.2}
Let $\tilde M = G\times_{G_Q, \rho} F$ be one 
of the manifolds given in Table 1 and adopt the same 
notation of \S 4.2; in particular, set  $N_F = N_F^{(1)}$ and 
$\varepsilon_F =\left\{\matrix 2 & \text{if}\ \  F = \CP^r\hfill\\
1 & \text{if}\ \  F = Q^r\hfill\endmatrix\right.$. 
\par
Consider a positive real  number 
$\hat c >0$, denote by $V_{\hat c}
 = \left(\frac{\frac{N_F}{\varepsilon_F} + 1}{\hat c}\right)^2$ and let
$V: [0, 1] \to [0, V_{\hat c}] $,
be a smooth map such that: 
\roster
\item"a)"
$\lim_{\theta \to 0} V(\theta) = 0$, 
$\lim_{\theta \to 0} \dot V(\theta)  = \dot V_o >0$ and $\lim_{\theta \to 1}
V(\theta) = 
V_{\hat c}$ (here  we set 
$\dot V = \frac{d V}{d\theta}$,
$\overset{\cdot \cdot}\to V = \frac{d^2 V}{d\theta^2}$
and   $\overset{\cdot \cdot \cdot}\to V = \frac{d^3 V}{d\theta^3}$);
\item"b)" the limits
$\lim_{\theta \to 1} \dot V(\theta) = \dot V_1 > 0$, $\lim_{\theta \to 1}
\overset{\cdot \cdot}\to V (\theta)$
and $\lim_{\theta \to 1} \overset{\cdot \cdot \cdot}\to V(\theta)$ 
exist and are finite;
\item"c)" $V$ verifies on $]0, 1[$ the equation
$$ \overset{\cdot \cdot}\to V
 + 
 \frac{ \dot V^2}{2\sqrt{V}}
\left(\frac{ N_F -1 }{\sqrt{V}} + 
\sum_{m = N_F+1}^{n-1} \frac{1}
{\sqrt{V} + a_m}\right) +  \dot V\left(\frac{ \hat c \sqrt{\frac{V}{\theta}}  
- \frac{N_F}{\varepsilon_F} - 1}{1-\theta} -
\frac{N_F - 1}{2\theta}\right) = 0
\tag 5.1$$
where $a_m = a_{\beta_m}$, with  $a_{\beta_m}$  
as defined in Theorem 4.2.
\endroster
Suppose also that the K-manifold $\tilde M$ 
is so that
\roster
\item"d)" for any $\beta_m \in R'_+$ the following inequalities hold: 
$$\left| \B(Z_D, Z_D) \cdot \frac{\B(Z^\kappa, i H_{\beta_m})}{\B(Z_\D, i
H_{\beta_m})}
\right| >  N_F + \varepsilon_F\ ,\qquad
\B(Z^\kappa, i H_{\beta_m}) < 0\ .$$
\endroster
Then  
the curve $Z_\omega(t) = f(t) H + I_\omega \in \z(\l) + \R Z_\D$, where $I_\omega
\in \z(l)$ is as in
(4.22)  and where 
$$f(t) =    \sqrt{V(\tanh^2(\varepsilon_F t))}\ ,\tag 5.2$$
is the algebraic representative 
of a K\"ahler-Einstein  metric  on $\tilde M$
with Einstein constant $c = -  \hat c\frac{\varepsilon_F}{\B(Z_\D, Z_\D)}$.\par
\endproclaim
\demo{Proof} Let $f(t) =  \sqrt{V(\tanh^2( \varepsilon_F t))}$ and 
$\theta(t) = \tanh^2(\varepsilon_F t)$. 
We claim that a) -  d) imply that 
$f$ verifies the hypothesis 1) - 4) of Theorem 4.2. \par
From definitions
$$\frac{d}{dt}\theta(t) = \frac{2 \varepsilon_F \tanh( \varepsilon_F 
t)}{\cosh^2(\varepsilon_F  t)} = 2 
\varepsilon_F \sqrt{\theta}(1 -
\theta) 
$$
and hence
$$f'(t) = \frac{\varepsilon_F}{\sqrt{V}} \dot V \sqrt{\theta}(1 - \theta)\ , \tag
5.3$$
$$f''(t) = - \varepsilon_F^2 \frac{(1 - \theta) ( - V \dot V + 
3 \theta V \dot V + \theta \dot V^2 - \theta^2 \dot V^2 - 
2 \theta V \overset{\cdot \cdot}\to V  + 
2 \theta^2 V \overset{\cdot \cdot}\to V )}{V^{3/2}}\ ,
\tag 5.4$$
$$f'''(t) = \varepsilon_F^3\left[- \frac{8 (1 -
 \theta)^2 \sqrt{\theta} \dot V}{\sqrt{V}} +
\frac{4 (1 - \theta) \theta^{3/2} \dot V}{\sqrt{V}} -
\frac{3 (1 - \theta)^3 \sqrt{\theta} \dot V^2}{V^{3/2}} +\right.$$
$$ +
\frac{6 (1 - \theta)^2 \theta^{3/2} \dot V^2}{V^{3/2}} +
\frac{3 (1 - \theta)^3 \theta^{3/2} \dot V^3} {V^{5/2}} +
\frac{6 (1 - \theta)^3 \sqrt{\theta} \overset{\cdot \cdot}\to V}{\sqrt{V}}  -
\frac{12 (1 - \theta)^2 \theta^{3/2} \overset{\cdot \cdot}\to V}{\sqrt{V}} - $$
$$\left. -
\frac{6 (1 - \theta)^3 \theta^{3/2} \dot V \overset{\cdot \cdot}\to V}
{\sqrt{V}^{3/2}} +
\frac{4 (1 - \theta)^3 \theta^{3/2} \overset{\cdot \cdot \cdot}\to
V}{\sqrt{V}}\right]\ .\tag 5.5$$
It follows that 
$$\frac{f''(t)}{f'(t)} = \varepsilon_F
\left[\left(2\frac{\overset{\cdot \cdot}\to V}{\dot V} -
\frac{\dot V}{V}
\right) (1-\theta) \sqrt{\theta} + 
\frac{1}{\sqrt{\theta}} (1 - 3 \theta)\right]\tag 5.6$$
$$f'(t)\left(\frac{ N_F }{f(t)} + 
\sum_{m = N_F+1}^{N} \frac{1}
{f(t) + a_m}\right) = \frac{\varepsilon_F}{\sqrt{V}} \dot V \sqrt{\theta}(1 -
\theta)
\left(\frac{ N_F }{\sqrt{V}} + 
\sum_{m = N_F+1}^{N} \frac{1}
{\sqrt{V} + a_m}\right)\tag 5.7$$
and 
$$\tilde c f(t)  - N_F  \tanh(t) - N_F 
\coth(t) = \tilde c f(t)  - 2 N_F  \coth(2t) = $$ 
$$ = 
 \tilde c \sqrt{V} - N_F (2 - \varepsilon_F) \sqrt{\theta} -  
\frac{\varepsilon_F N_F }{\sqrt{\theta}}\ .\tag 5.8$$
After some simple algebraic manipulation, 
it follows immediately that $V$ verifies (5.1) 
if and only if $f(t)$ verifies (2) of Theorem 4.2
with 
$\tilde c = 2 \varepsilon_F \hat c$. \par
Moreover, if a) and b) are verified, from (5.6) - (5.8) we get that 
$$\lim_{t\to 0} f(t) = \lim_{\theta\to 0} \sqrt{V(\theta)} = 0$$
$$\lim_{t\to 0} f'(t) = \dot V_o \sqrt{\lim_{\theta\to 0}\frac{\theta}{V}}
= \sqrt{\dot V_o} >0$$
Similarly, one can check that $\lim_{t\to 0} f''(t) = 0$ and that
$\lim_{t\to 0} f'''(t)$ is finite, so that all parts of condition (3) of Theorem 4.2 is
verified, with the only exception of the inequalities $\B(Z^\kappa, i H_{\beta_m}) < 0$.
\par
By some tedious but direct computations, one can check in the same way 
that also all parts of condition (4) of Theorem 4.2 are verified, with
the only exception of the positivity of the limits 
$\lim_{t\to + \infty} f(t)$ and 
$\lim_{t\to + \infty} - (f(t) + a_m) \beta_m(Z_\D)$.
\par
\bigskip
So, in order to conclude the proof, we need to check that if
d) is true and  $V$ verifies the 
boundary conditions a), then the condition (1) and 
the remaining parts of (3) and 
(4) in Theorem 4.2 are automatically 
verified by $f(t) = \sqrt{V(\tanh^2(t))}$. \par
For this we need a pair of preparatory lemmata.\par
\enddemo
\bigskip
Let us introduce the notation 
$\alpha = \frac{N_F -
1}{2}\geq 0$ and $\alpha' = \frac{\frac{N_F}{\varepsilon_F} -
1}{2}\geq  - \frac{1}{4}$; note  that $V_{\hat c} = 4\left(\frac{\alpha'+1}{\hat
c}\right)^2 > 0$.\par
Let us 
also rewrite 
(5.1) as
$$\frac{d}{d\theta}\left[\log\left(\dot V V^\alpha \prod_{m =
N_F+1}^{n-1}|\sqrt{V} + a_m|\right)\right] = 
\frac{2(1 + \alpha') - \hat c \sqrt{\frac{V}{\theta}}}{1-\theta} +
\frac{\alpha}{\theta}\tag 5.9$$
\medskip
\remark {Remark 5.3} Note that, since $\hat c = - c \frac{\B(Z_\D,
Z_\D)}{\varepsilon_F}$,
if we assume that for $0 < \theta < 1$
$$0 < V(\theta) <
V_{\hat c} = 4 \left(\frac{\alpha'+1}{\hat c}\right)^2$$
and that  hypothesis d) holds true,  
 then the sign of 
$a_m$ and $a_m + \sqrt{V(\theta)}$ does not depend on $\theta$
(to check this, use Table A1 in Appendix);  in
particular, 
the value of $a_m + \sqrt{V(\theta)}$  is never equal to $0$.\endremark
\medskip
Now, from (5.9) and the above remark, we immediately obtain the
following lemma.
\bigskip
\proclaim{Lemma 5.4} Assume that d) of Proposition 5.2 is true
and let $V: [a,b] \subset ]0,1[ \to \R$ be a $\Cal C^1$-solution 
of (5.1) with 
$ 0 
< V(\theta) < V_{\hat c} = 4\left(\frac{\alpha'+1}{\hat c}\right)^2$
for any $\theta \in [a,b]$. 
Then, for $\theta_1 \leq \theta_2 \in [a,b]$,
$$\left(\frac{\theta_2}{\theta_1}\right)^\alpha
\left(\frac{1 + \sqrt{\theta_1}}{1 + \sqrt{\theta_2}}\right)^{4(1 + \alpha')}
< \frac{\dot V({\theta_2})}{\dot V({\theta_1})}
\left(\frac{V({\theta_2})}{V({\theta_1})}\right)^{\alpha}
\prod_{m = N_F +1}^{n-1}\left(\frac{|a_m + \sqrt{V({\theta_2})}|}{|a_m +
\sqrt{V({\theta_1})}|}\right)
< $$
$$ <\left(\frac{\theta_2}{\theta_1}\right)^\alpha\left(\frac{1-\theta_1}{1 -
\theta_2}\right)^{2(1+\alpha')}\tag5.10$$
\endproclaim
\demo{Proof} By hypothesis, 
$$2(1 + \alpha')\frac{\left(1 -
\frac{1}{\sqrt{\theta}}\right)}{1-\theta}  <
\frac{2(1 + \alpha') - \hat c\sqrt{\frac{V}{\theta}}}
{1-\theta} < \frac{2(1+ \alpha')}{1-\theta}$$
From these inequalities and from (5.9), the claim follows immediately by
integration.\qed
\enddemo
\bigskip
\proclaim{Lemma 5.5}  Assume that hypothesis d) holds true and 
let $V: [a,b] \subset ]0,1[ \to \R$ be a $\Cal C^1$-solution 
of (5.1) with $0  < V(a) < V(b) <4\left(\frac{1 + \alpha'}{\hat c}\right)^2$.
Then
$V$ is monotone 
increasing on  $[a,b]$ and 
$0 < V(\theta) < 4\left(\frac{1 + \alpha'}{\hat c}\right)^2$ for any $\theta\in
[a,b]$.
\endproclaim
\demo{Proof} By hypothesis, there exists at least one point $\theta_o \in ]a, b[$
such that $0 <V({\theta_o}) < V(b) < 4\left(\frac{1 + \alpha'}{\hat c}\right)^2$ 
and  with $\dot V({\theta_o}) >0$.
Moreover, by Remark 5.3, 
 we also have that  $\prod_m\left( a_m + \sqrt{V({\theta_o})}\right) \neq 0$.\par
Now, consider the set $K = \{\ \theta \in [\theta_o, b]\ :\ \dot V(\theta) >0\ ,
\ V(\theta) < 4\left(\frac{1 + \alpha'}{\hat c}\right)^2\ \}$. It is clear that $K$
is 
open in $[\theta_o, b]$. Moreover, if $\theta_1 = \sup K$, 
 from the left hand side of (5.10) it follows that
$$M_{\theta_o}\left(\frac{1 + \sqrt{\theta_o}}{1 +
\sqrt{\theta_1}}\right)^{4(\alpha'+1)} \left(\frac{\theta_1}{\theta_o
}\right)^\alpha
<  {\dot V}(\theta_1)
V(\theta_1)^{\alpha}
\prod_m |a_m + \sqrt{V(\theta_1)}|
\tag 5.11$$
with 
$M_{\theta_o} = {\dot V}(\theta_o)
V(\theta_o)^\alpha
\prod_m |a_m + \sqrt{V(\theta_o)}|\neq 0$. This shows that
${\dot V}(\theta_1) > 0$ and hence that it belongs to $K$.
In particular, it follows that $K$ is open and closed and hence it coincides with
$[\theta_o, b]$, proving the claim for any $\theta \in [\theta_o, b]$.\par 
A similar argument, which uses the right
hand side of (5.10) in place of the left hand side, brings to the same conclusion
for any $\theta \in [a,\theta_o]$.\qed
\enddemo
\bigskip
\demo{End of Proof of Proposition 5.2} To conclude the proof, it suffices to
observe that if
$V: [0, 1] \to [0, 
V_{\hat c}]$, 
with $V_{\hat c} = 4(\alpha' + 1)^2/\hat c^2$, then Lemma 5.5 applies on
$V|_{[a,b]}$ for any 
interval $[a,b] \subset ]0,1[$. In particular, $V$ is strictly increasing and 
$f(t) = V(\tanh^2(t))$ is monotone increasing and  positive on $]0, \infty[$. So
it remains 
to check that the positivity of 
$- (f(t) + a_m) \beta_m(iZ_\D)$ for any $t \in ]0, \infty[$ and of its 
limits for $t \to 0$ and $t\to + \infty$;
this 
is an immediate consequence of the hypothesis d) and Remark 5.3.\qed
\enddemo
\bigskip
The following theorem shows that, if d) of Proposition 5.2 holds true, then a
solution 
for the K\"ahler-Einstein equation  exists with 
only one possible exception. The proof of this fact is
practically an 
adaptation to our more general case
 of the arguments used in \S 6 and \S 7 of \cite{GC}
\footnote{We are indebted to D. Guan for pointing us a serious 
mistake in a previous version of the proof of Theorem 5.6}.\par
\medskip 
\proclaim{Theorem 5.6} Under the hypothesis d) of Proposition 5.2 and 
assuming that
 $\tilde M \neq \SO_{10}\times_{(\operatorname{T}^1\times \SO_8)} Q^7$, 
for any constant $\hat c >0$, 
there exists a real analytic solution
 $V: [0, 1] \to [0, V_{\hat c}]$ of the differential problem given 
by conditions a), b) and c) of Proposition  5.2 
\endproclaim
\demo{Proof} The proof is subdivided into several claims.\par
\medskip
\noindent{\bf Claim 1}. \  {\it For any interval $[a,b] \subset ]0, 1[$ and
any
positive constant  $0 < c_a   < V_{\hat c}$ there exists 
a monotone increasing solution $V: [a, b] \to ]0, V_{\hat c}[$  of (5.1) with
$V_a
= c_a$\/}.
\par
\smallskip
\medskip
In order to prove this claim, we consider the Cauchy problem ($C_{k}$)
given by the differential equation (5.1) together with initial conditions
$V(a) = c_a$ and ${\dot V}(a) = k > 0$. Let us denote by $[a,b_k[\subset [0,1[$
the maximal interval on which there exists a solution $V_k$ of ($C_k$),
satisfying 
$V_k(\theta) < V_{\hat c}$ for all $\theta \in [a,b_k[$. By Lemma 5.5, the
solution $V_k$ is monotone increasing. We also have that $\lim_{\theta\to
b_k}V_k(\theta)
= V_{\hat c}$, otherwise the interval $[a,b_k[$ would not be maximal. \par
We claim that
$\sup_{k>0}b_k = 1$. Indeed, suppose that $\sup_{k>0}b_k = b^* < 1$. It would
imply that
for 
every $k>0$, there exists some $\theta_k \in ]a,b_k[$ so that 
$${\dot V}_k(\theta_k) = {{V_{\hat c} - c_a}\over {b_k - a}} >
{{V_{\hat c} - c_a}\over {b^* - a}}.$$
However, from Lemma 5.4, it would follow that there exists some 
positive constant $H$ so that
$$0 < {{V_{\hat c} - c_a}\over {b^* - a}} < k H \left({{1-a}\over{1-b^*}}\right)
^{2(1+\alpha')} \left({{b^*}\over{a}}\right)^\alpha\tag 5.12$$
for any $k$. This leads to a contradiction when we let $k \to 0$.
\medskip
\noindent{\bf Claim 2}. \  {\it  Let $V: [a,b] \subset ]0, 1[ \to
]0, V_{\hat c}[$
be a monotone increasing solution of (5.1) with $V(a) >0$. 
Then for any positive number $c_b$ with $V(b) < c_b < V_{\hat c}$,
there exists a monotone increasing solution $\tilde V$ of (5.1) with 
$\tilde V(a) = V(a)$ and $\tilde V(b) > c_b$\/}.\par
\smallskip
As before, let us denote by ($C_{k}$) the Cauchy problem 
given by the differential equation (5.1) together with initial conditions
$V(a) = c_a$ and ${\dot V}(a) = k > 0$.
Let also $K\subset [0,+\infty[$ be the set of all $k \geq 0$ such that
the Cauchy problem ($C_k$) has a solution $V_k$ defined on $[a,b]$ and with 
$V_k(b) < V_{\hat c}$.\par
The claim is proved if we can show that $\sup_{k\in K}V_k(b) =
V_{\hat c}$. Suppose not and let  $\sup_{k\in K}V_k(b) =
\lambda < V_{\hat c}$. In particular, we have that for any $k\in K$, any solution
$V_k$ 
of ($C_k$) on $[a,b]$ verifies $V_k(b)\leq \lambda < V_{\hat c}$ and hence, using
(5.10), 
one can easily check that  $K$ is open and closed
in $[0,+\infty[$ and  hence that $K = [0,+\infty[$. But then,  using the same
arguments
as in the previous claim, we obtain that 
$$k \left({{b}\over{a}}\right)^{\alpha}
\left({{1+\sqrt{a}}\over{1+\sqrt{b}}}\right)^{4(1+\alpha')}
\left({{c_a}\over \lambda}\right)^\alpha \prod_{m=N_F+1}^{n-1}
\left|{{a_m+\sqrt{c_a}}\over{a_m+\sqrt{\lambda}}}\right| \leq {{\lambda - c_a}\over{b - a}}\
,\tag 5.13$$
for any $k>0$, which is impossible. \par
\bigskip
\noindent{\bf Claim 3}. \  {\it For any $k \in ]0,1[$ there exists a monotone 
increasing real analytic
solution  $V: ]0,k]  \to ]0, V_{\hat c}]$ of (5.1) with \/}
$$\lim_{\theta \to 0} V(\theta) = 0\ ,\qquad \lim_{\theta \to k} V(\theta) =
V_{\hat c}$$
\par
\smallskip
Using Claim 1 and 2, we may construct a sequence $\{V^{(n)}\}$ of monotone
increasing solutions
of (5.1) defined on the intervals $[\frac{1}{n}, k - \frac{1}{n}]$, with values
in 
$]0, V_{\hat c}[$  and such that 
$$
\lim_{n \to \infty} V^{(n)}(\frac{1}{n}) = 
\lim_{n \to \infty}\frac{1}{n} = 0\ ,\quad \lim_{n\to \infty} V^{(n)}(k -
\frac{1}{n}) = V_{\hat
c}\ .$$
Note that, for any interval $[a, b] \subset ]0, k[$, there exists an integer
$N_{a,b}$
such that 
for any $n > N_{a,b}$, 
$$[a,b] \subset [\frac{1}{n}, k - \frac{1}{n}]\ .$$
Now, fix a positive value 
$0< \varepsilon < a$ and notice that, for any $\theta \in [a,b]$ and any function
$V^{(n)}$, 
with $n$ large enough that $[a - \varepsilon,b] \subset [\frac{1}{n}, k -
\frac{1}{n}]$, 
there exists a
value $a - \varepsilon < \theta_n <  \theta$ so that 
$$\dot V^{(n)}(\theta_n) = \frac{V^{(n)}(\theta) - V^{(n)}(a - \varepsilon)
}{\theta - a +  \varepsilon}
< \frac{V^{(n)}(b) - V^{(n)}(a - \varepsilon) }{\varepsilon} < \frac{V_{\hat
c}}{\varepsilon}\ .$$
From  (5.10) (where we replace $\theta_1$ and $\theta_2$ with 
$\theta_n$ and $\theta$, respectively) and from the fact that for any $n$, 
$\frac{V^{(n)}(\theta_n)}{V^{(n)}(\theta)} < 1$, it follows 
that there exists a constant $C$, independent of $n$ and the interval 
$[a,b]\subset ]0,k[$, so that
$\dot V^{(n)}(\theta) < C$
for any $n$ sufficiently large and 
for any $\theta \in [a,b]$. \par
In particular, we conclude that the functions
$\dot V^{(n)}|_{[a,b]}$
are uniformly bounded  and hence that the functions $V^{(n)}|_{[a,b]}$ converge,
up to a subsequence, to a $\Cal
C^0$-function, say $V$.\par
Consider now an exhausting sequence of intervals $[a_n,b_n] \subset ]0, k[$: on
each of these intervals
we may define the function $V: [a_n, b_n]\to ] 0, V_{\hat c}[$ as limit of the
sequence $V^{(n)}$; 
since the limit function $V$ coincides on the intersection of two intervals
$[a_{n_1}, b_{n_1}], 
[a_{n_2}, b_{n_2}]$, we obtain a nondecreasing continuos function $V$ defined 
on the open
interval $]0, k[$. Moreover, since 
$$\lim_{\theta\to 0}V(\theta) 
= \inf_{\theta\in ]0,k[}V(\theta) 
\leq \inf_{\theta\in ]{1\over n},k-{1\over n}]}V^{(n)}(\theta) 
= {1\over n},$$
we get $\lim_{\theta\to 0}V(\theta)  = 0$. Now fix $\varepsilon > 0$ and 
$t_o > k - {\varepsilon\over{3C}}$; take also $N_o$ large enough so that $t_o \leq k - {1\over n}$
and $|V^{(n)}(k-{1\over n}) - V_{\hat c}| < \varepsilon/3$ for all $n \geq N_o$. Then,
choosing a suitable large $n$
$$\eqalign{|V(t_o) - V_{\hat c}| \leq & |V(t_o) - V^{(n)}(t_o)| + |V^{(n)}(t_o) - 
V^{(n)}(k-{1\over n})| +\cr
{}& |V^{(n)}(k-{1\over n}) - V_{\hat c}| \leq |V(t_o) - V^{(n)}(t_o)| + C|t_o - k + 
{1\over n}| + {\varepsilon\over 3} \leq \varepsilon}$$
Since  $V$ is non decreasing,  from this we have that 
$\lim_{\theta\to k}V(\theta) = V_{\hat c}$.\par
To conclude the proof of the claim, it remains to check that 
$V$ is actually a solution of (5.1), and then it will also follow that 
it is  real analytic. \par
\medskip
First of all, notice that for any $\theta \in ]0, k[$, the value $V(\theta) \neq
0$.\par 
In fact, assume that there is some $\theta_o$ so that $
\lim_{n\to 0} V^{(n)}(\theta_o) = V(\theta_o) = 0$. On the other hand, since
$\lim_{\theta \to k} V(\theta) >0$, 
there 
has to be  some point $\theta_o + \delta \in [\theta_o, k[$ so that $V(\theta_o +
\delta) > 0$. So, 
we may consider the functions $V^{(n)}|_{[\theta_o - \delta', \theta_o +
\delta]}$ for some 
fixed $0 < \delta' < \theta_o$. 
Fix $\varepsilon >0$; since  each function $V^{(n)}$ is positive and strictly
increasing,  we have that 
for all $n$ sufficiently large, there exists some $\theta_n \in ]\theta_o -
\delta', \theta_o[$ and 
some $\xi_n \in ]\theta_o , \theta_o + \delta[$ such that 
$$V^{(n)}(\theta_o) < \varepsilon\ ,$$
$$\dot V^{(n)}(\theta_n) V^{(n)}(\theta_n)^\alpha = \frac{V^{(n)}(\theta_o) -
V^{(n)}(\theta_o - \delta')}{\delta'}
V^{(n)}(\theta_n)^\alpha
< \frac{V^{(n)}(\theta_o)^{1 + \alpha}}{\delta'}\ , $$
$$\dot V^{(n)}(\xi_n) V^{(n)}(\xi_n)^\alpha =  \frac{V^{(n)}(\theta_o + \delta) -
V^{(n)}(\theta_o)}{\delta} V^{(n)}(\xi_n)^\alpha 
>$$
$$ > \frac{V(\theta_o + \delta)  -  \varepsilon}{\delta}
V^{(n)}(\theta_o)^\alpha\ . $$
From these inequalities, we get that for any $\varepsilon >0$ sufficiently small,
 there exists some 
$\theta_o - \delta' < \theta_n < \xi_n 
< \theta_o + \delta$ such that 
$$\frac{\dot V^{(n)}(\xi_n) V^{(n)}(\xi_n)^\alpha}{
\dot V^{(n)}(\theta_n) V^{(n)}(\theta_n)^\alpha} >
\delta' \frac{V(\theta_o + \delta)  -  \varepsilon}{\delta V^{(n)}(\theta_o)}  >
\frac{1}{\varepsilon} 
\frac{\delta'}{\delta} \frac{V(\theta_o + \delta)}{2}$$
which brings to an immediate contradiction with the inequalities (5.10).\par
\medskip
At this point, it remains to observe that, being $V(\theta) = \lim_{n \to \infty}
V^{(n)}(\theta) 
\neq 0$ for any $\theta \in ]0,k[$ and 
since, for any closed interval $[a, b] \subset ]0,k[$, the first derivatives of
the 
functions $V^{(n)}|_{[a,b]}$ are uniformly bounded, it follows from (5.1) that
also 
the second derivatives  of
 $V^{(n)}|_{[a,b]}$ are uniformly bounded. So, for any interval $[a,b]$, 
the sequence
$V^{(n)}|_{[a,b]}$ converges to $V|_{[a,b]}$  in  $\Cal C^1$ and we get that 
the limit function $V$ is a solution 
of  (5.1), by  well known results on the
smooth dependence of solutions from the  initial data. 
\par
\bigskip
\noindent{\bf Claim 4}. \  {\it In case
$\tilde M \neq \SO_{10}\times_{(\SO_2\times \SO_8)} Q^7,\ 
\operatorname{E}_6 \times_{(\SO_2\times \operatorname{Spin}_{10})} Q^9$,   
there exists a monotone 
increasing real analytic
solution  $V: ]0,1[  \to ]0, V_{\hat c}[$ of (5.1) with \/}
$$\lim_{\theta \to 0} V(\theta) = 0\ ,\qquad \lim_{\theta \to 1} V(\theta) =
V_{\hat c}\ .$$
\par
For any $k \in \Bbb N$ denote by $V_k : [0,1-{1\over k}] \to [0,V_{\hat  c}]$ the solution 
of (5.1) given by Claim 3. Note that the same arguments used in the 
proof of Claim 3 show that  the values $|\dot V_k(\theta)|$ are uniformly bounded 
on intervals $]0,b]$ for $b <1$. This means that the sequence 
$\{V_k\}$ converges uniformly on compacta  of $[0,1[$ to a nondecreasing solution of 
(5.1) 
with $V:]0,1[ \to ]0, \hat V_{\hat c}[$ and that 
$\lim_{\theta\to 0} V(\theta) = 0$. It remains to check that 
$\lim_{\theta\to 1} V(\theta) = \sup_{\theta\in [0,1[}V(\theta)$ is equal to 
$V_{\hat c}$.\par
We fix $0 <\varepsilon < \min\{(\alpha' + 1), 1/3\}$ and for any $k$ we choose $\theta_k \in ]0,1-{1\over k}[$ 
with $V_k(\theta_k) = \left({{2\alpha' + 2 - 2\varepsilon}\over {\hat c}}\right)^2$. 
Note that $0$ is not a limit point for $\{\theta_k\}$: this follows from the mean 
value theorem together with the fact that $\dot V_k$ are uniformly bounded on 
the compact subsets of $]0,1/2[$ (see proof of Claim 3). \par
\bigskip
\proclaim{Lemma 5.7} In case
$\tilde M \neq \SO_{10}\times_{(\operatorname{T}^1\times \SO_8)} Q^7$, 
for any sufficiently small $\varepsilon > 0$, we have that
 $\sup_k \{\theta_k\}<1$.
\endproclaim
\demo{Proof} We consider the change of variable $\theta(s) = 1 - e^{-s}$.
Since
$\frac{d\theta}{ds} = e^{-s} = 1 - \theta$, for any function $\Cal F_\theta$ we
have that 
$$\Cal F'_s \overset\operatorname{def}\to= \left.\frac{d{\Cal F}}{ds}\right|_{s}
= \dot {\Cal
F}(\theta(s)) (1 - \theta(s))\ .$$
From (5.9), we  obtain that 
$$\frac{d}{ds}\left[\log\left(V_k' V_k^\alpha \prod_{m = N_F+1}^{n-1}(|\sqrt{V_k} + a_m|)
\right)\right] = \alpha \left({{1-\theta(s)}\over{\theta(s)}}\right) + 1 + 2 \alpha'  - \hat
c\sqrt{\frac{V_k}{\theta(s)}} 
\tag 5.14$$
We consider the points $s_k\in \R$ such that $\theta(s_k) = \theta_k$ and, up to 
a subsequence, we suppose 
that the sequence of points $\theta(s_k)$ tends to $1$, that is 
$\lim_{k\to +\infty}s_k = +\infty$, aiming to obtain a contradiction. \par
We consider the functions $\tilde V_k(s) \= V_k(\theta(s + s_k))$. First of all, 
we want to prove that, for every compact $K \subset \R$,  the derivatives $\tilde V_k'|_K$ 
are uniformly bounded for all $k$ large enough. Observe that the functions $\tilde V_k$ satisfy 
the equation
$$\eqalign{
\tilde V_k(s)'' &= (\tilde V_k'(s))^2 \left(\sum_l{1\over{\al^2 - \tilde V_k(s)}} - 
{{\alpha}\over {\tilde V_k(s)}}\right) + \tilde V_k'(s)
\left(\left.{{\alpha(1 - \theta)}\over\theta}\right|_{s+s_k}\right. +\cr
{}& + \left.1 + 2\alpha' - 
\hat c {{\sqrt{\tilde V_k(s)}}\over {\sqrt{\theta(s+s_k)}}}\right)\ ,} \tag 5.15$$
where we denote by $\al$ the positive coefficients which appear in the 
collection of coefficients $\{a_\ell\}$; indeed, 
we recall that the coefficients $a_l$ appear in pairs of 
opposite signs (see Table A1 in Appendix).\par
Note that $\theta(s + s_k) \geq \theta(s_k)$ for $s\geq 0$;  hence, for $k$ large enough, 
we have that for all $s\geq 0$ 
$$\tilde V_k''(s) \leq C (\tilde V_k'(s))^2 + \tilde V_k'(s) (-1 + 3\varepsilon)\ ,$$
where $C = \sum_l{1\over{\al^2 - V_{\hat c}}}$. By integration, since $\tilde V_k(0)$ 
does not depend on $k$, we have for all $s \geq 0$
$$\tilde V_k'(s) \leq C_1 \tilde V_k'(0) e^{(3\varepsilon-1)s}\tag 5.16$$
for some positive constant $C_1$. \par
We claim now that the sequence of values $\{V_k'(0)\}$ is bounded from above and 
from below by two positive constants. \par
Indeed, integrating (5.16), we get that
$$\tilde V_k(s) - \tilde V_k(0) \leq C_1 \tilde V_k'(0)\int_0^\infty 
e^{(3\varepsilon-1)t}dt$$
and therefore $V_{\hat c} - \left({{2\alpha' + 2 - 2\varepsilon}\over {\hat c}}\right)^2
\leq {{C_1}\over{1-3\varepsilon}} \tilde V_k'(0)$, showing that 
$\inf_k \tilde V_k'(0)$ is positive. We now show that $\tilde V_k'(0)$ are 
also bounded from above. Note that 
$$\left.{{\alpha(1 - \theta)}\over\theta}\right|_{s+s_k} + 1 + 2\alpha' - 
\hat c { {\sqrt{\tilde V_k(s)}}\over {\sqrt{\theta(s+s_k)}} } \geq 
1 + 2\alpha' - \hat c {{\sqrt{\tilde V_{\hat c}}}\over {\sqrt{\theta(s+s_k)}}} 
\geq -{3\over 2}$$
for all $s\geq 0$ and all large enough $k$. Hence 
$$\tilde V_k''(s) \geq -{\alpha\over{\tilde V_k(s)}}\tilde V_k'(s)^2 - {3\over 2}
\tilde V_k'(s)$$ 
and by integration 
$$\tilde V_k'(s) \tilde V_k(s)^\alpha \geq C_2 \tilde V_k'(0) e^{-3s/2}\ ,$$
where $C_2 = \tilde V_k(0)^\alpha$ does not depend on $k$. Again, by integration, 
for all $s \geq 0$
$${1\over{\alpha+1}} V_{\hat c}^{\alpha+1} \geq 
{1\over{\alpha+1}}(\tilde V_k(s)^{\alpha+1} - \tilde V_k(0)^{\alpha+1}) \geq 
{{2C_2}\over 3} \tilde V_k'(0) (1 - e^{-3s/2})$$
and this shows that the sequence $\{\tilde V_k'(0)\}$ is also  bounded from above. \par 
\medskip
Let us now fix a compact interval $I = [-A^2,0]$. Since $s \geq - A^2$ 
for any $s\in I$, 
we have $\theta(s+s_k) \geq \theta(s_k - A^2) \to 1$
for $k\to +\infty$. Hence for $k$ large enough and for all $s \in I$, the right handside
of (5.14) is bigger or equal than $-2$ and, by integration,
$${{\tilde V_k'(0)(\prod_l \al^2 - \tilde V_k(0))\tilde V_k(0)^\alpha}
\over {\tilde V_k'(s)(\prod_l \al^2 - \tilde V_k(s))\tilde V_k(s)^\alpha}} 
\geq e^{-2\int_s^0 dt} = e^{2s}\ .$$
From this we conclude that, for some positive constants $C_3,C_4$,  
$$C_3 \tilde V_k'(0) e^{-2s} \geq 
\tilde V_k'(s)(\sum_l \al^2 - \tilde V_k(s))\tilde V_k(s)^\alpha \geq 
C_4 \tilde V_k'(s)\tilde V_k(s)^\alpha\ . \tag 5.17$$
From (5.17), 
it follows that  the functions $(\tilde V_k(s)^{\alpha+1})'$ are uniformly 
bounded on every compact subset in $]-\infty, 0[$. Now, from (5.16), the fact that 
the sequence $\{\tilde V_k'(0)\}$
is bounded from above and that the functions
$\tilde V_k(s)$ are uniformily bounded, it follows also that
the functions $(\tilde V_k(s)^{\alpha+1})'$ are uniformly 
bounded   
 on $[0,+\infty[$.  We conclude that, up to a 
subsequence, $\tilde V_k(s)^{\alpha+1}$ converge to a continuous, nondecreasing function $W:\R 
\to [0,V_{\hat c}]$. Note that $W(0) \neq 0$ and that, being 
$\tilde V_k'(0)$ bounded 
from below by a positive constant, the function $W$ is not constant.\par
\medskip
 We now prove that 
$W$ never vanishes; indeed, for any $s \in [-A^2,0]$, the right hand side of (5.14) does 
not exceed a constant $M$ and therefore 
$${{\tilde V_k'(0)(\prod_l \al^2 - \tilde V_k(0))\tilde V_k(0)^\alpha}
\over {\tilde V_k'(s)(\prod_l \al^2 - \tilde V_k(s))\tilde V_k(s)^\alpha}} \leq 
e^{Ms}\ .$$
Hence there are constants $C_5$, $C_6$ such that for all $s \in I$,
$$C_5 \tilde V_k'(0) \leq e^{Ms} \tilde V_k'(s) \tilde V_k(s)^{\alpha} \prod_l(\al^2 - 
\tilde V_k(s)) \leq C_6 \tilde V_k'(s) \tilde V_k(s)^{\alpha}.$$
Since $\inf_k\tilde V_k'(0) > 0$, we see that $(\tilde V_k(s)^{\alpha+1})'$ are 
bounded from below by a positive constant which  depends only on the interval $I$.
In particular $W'$ is never 
zero. On the other hand, if there exists some
$t_o < 0$ such that $W(t_o) = 0$, then $W(s) = 0$ for all $s < t$ and this is not 
possible. So, $W$ never vanishes and the sequence $\tilde V_k$ converges uniformly on compacta 
to a differentiable function $\tilde V \= (W)^{1/(\alpha+1)}$. Note that the 
limit function $\tilde V$ is not 
constant and nondecreasing. \par
\medskip
From (5.16) and the bounds on the sequence $\{\tilde V'_k(0)\}$, we see that 
$\tilde V'(s) \leq C_7 e^{(3\varepsilon - 1)s}$ for some constant $C_7$ and for all $s\geq 0$. 
Therefore $\lim_{s\to +\infty} \tilde V'(s) = 0$. \par
The function $\tilde V$ satisfies the 
equation 
$$\left( \log (\tilde V' \tilde V^\alpha \prod (\al^2 - \tilde V)\right)' = 
1 + 2\alpha' - \hat c \sqrt{\tilde V}\ . \tag 5.18$$
If we denote by $U(s) \= \tilde V' \tilde V^\alpha \prod (\al^2 - \tilde V)$, the equation (5.18) becomes
$$U' = (1 + 2\alpha' - \hat c\sqrt{\tilde V}) \tilde V^\alpha 
\prod (\al^2 - \tilde V) \tilde V'.\tag 5.19$$
Let us consider the function 
$$g(x) = (1 + 2\alpha' - \hat c\sqrt{x}) x^\alpha 
\prod (\al^2 - x)$$
 for $x \in [0,+\infty[$ and let $G(x) = \int_0^x g(t) dt$ for 
$x\geq 0$. Integrating (5.19) we get that 
for any  $s_1 , s_2 \in \R$
$$U(s_2) - U(s_1) = G(\tilde V(s_2)) - G(\tilde V(s_1))\ .\tag 5.20$$ 
Let us denote by 
$$a = \lim_{s\to -\infty}\tilde V(s)\ ,\qquad b = \lim_{s\to +\infty}
\tilde V(s)\ ,$$
and note that $a,b \in [0,V_{\hat c}]$.  \par
We already noticed that  $\lim_{s\to + \infty}
\tilde V'(s) = 0$ and therefore $\lim_{s\to + \infty} U(s) = 0$. Hence, 
from (5.20) it follows that $- U(s) = G(b) - G(\tilde V(s))$. This implies also that 
 $\lim_{s\to -\infty}U(s)$ exists and it is finite. Indeed, we claim that $\lim_{s\to -\infty}U(s) = 0$:
in fact,  in case $\lim_{s\to -\infty}U(s)  = \lambda \neq 0$, from the existence of 
$\lim_{s\to - \infty}\tilde V(s)$ and the definition of $U$, we  infer that also 
$\lim_{s\to - \infty}\tilde V'(s)$ exists; but this implies that 
$\lim_{s\to - \infty}\tilde V'(s) = 0$ 
and it contradicts the hypothesis that $\lambda \neq 0$.
\medskip
From (5.20) and the previous remarks, we conclude that
$$0 = G(b) - G(a) = \int_a^b g(x) dx \ .\tag 5.21$$
Now, if we can prove that $\int_0^{V_{\hat c}} g(x) dx < 0$, 
 we immediately get a contradiction and we conclude the proof of the lemma.  Indeed, 
notice that if we choose $\varepsilon$ small enough,
we have that the value
$$b = \lim_{s \to + \infty} \tilde V(s) > \tilde V(0) = V_k(\theta_k) = 
4 \left(\frac{\alpha'+1 - \varepsilon}{\hat c}\right)^2$$
can be made arbitrarily close to $V_{\hat c} = 4 \left(\frac{\alpha'+1}{\hat c}\right)^2$;
in particular, we can choose $\varepsilon$ so that $\int_0^b g(x) dx < 0$. Now, 
since $g(x) \geq 0$ if and only if $x \leq \left({{1+2\alpha'}\over{\hat c}}\right)^2$, 
we have that in case $a \geq \left({{1+2\alpha'}\over{\hat c}}\right)^2$, the integral 
$\int_a^b g(x) dx$ is negative, which is contradictory with (5.21); when 
$a < \left({{1+2\alpha'}\over{\hat c}}\right)^2$,we have that
$\int_a^b g(x) dx < \int_0^b g(x) dx <0$ and again this contradicts (5.21). 
\par
\medskip
In the Appendix  we estimate the sign of the integral $\int_0^{
V_{\hat c}}g(x) dx$  for all cases of Table 1. The reader can check that 
for all possibilities, except the cases 
$\tilde M =  \SO_{10}\times_{(\operatorname{T}^1\times \SO_8)} Q^7,\ 
\operatorname{E}_6\times_{(\SO_2\times \operatorname{Spin}_{10})} Q^9$, 
this integral is negative. This concludes the proof. \qed
\enddemo
\bigskip
We may now conclude the proof of Claim 4. 
For a given small $\varepsilon > 0$, by Lemma 5.7, we may suppose that 
$\lim_{k\to \infty} \theta_k = \theta_o < 1$. Therefore for $\theta > \theta_o$, we 
have that $V_k(\theta)$ is bigger than 
$\left({{2\alpha' + 2 - 2\varepsilon}\over {\hat c}}\right)^2 = V_k(\theta_k)$ for $k$ 
large enough. So $\lim_{\theta\to 1} V(\theta) = \sup_{\theta\in [0,1[}V(\theta)
\geq \left({{2\alpha' + 2 - 2\varepsilon}\over {\hat c}}\right)^2$. 
By the freedom on the choice 
of $\varepsilon >0$, we have $\lim_{\theta\to 1} V(\theta) = V_{\hat c}$.
\bigskip
\noindent{\bf Claim 5}. \  {\it Let $V: ]0,1[ \to ] 0, V_{\hat c}[ $
be the solution of (5.1) given in Claim 4. Then, for any 
 $\varepsilon>0$ there exist a positive constants $M_{\varepsilon}$
and a point 
$\theta_\varepsilon \in ]0,1[$ so that,  for any 
$\theta \in
[\theta_\varepsilon, 1[$ \/}
$$\dot V(\theta) (1- \theta) < 
M_\varepsilon
(1- \theta)^{1 - 2 \varepsilon}\
.\tag 5.22$$
\par
\smallskip
We use the same change of variable
$\theta(s) = 1 - e^{-s}$ we used in the proof of Lemma 5.7. 
Since $\lim_{\theta \to 1} \hat c\sqrt{\frac{V}{\theta}} = \hat c \sqrt{V_{\hat
c}} = 2(\alpha' + 1)$, 
 it follows that for any given $\varepsilon >0$, there exists some $s_o$ such
that 
for all $s \geq s_o$ (and hence for all $\theta(s) \geq \theta_\varepsilon =
\theta(s_o)$), we have that 
$$
\frac{d}{ds}\left[\log\left(V' V^\alpha \prod_{m = N_F+1}^{n-1}|\sqrt{V} + a_m|
\left(\frac{1}{\theta(s)}\right)^\alpha\right)\right]  < 1 + 2 \alpha' -
2(\alpha' + 1) + 2 \varepsilon = 
-1 + 2\varepsilon \ . $$
By integration, we get that for any $\theta_{\varepsilon} \leq  \theta < 1 $
$$\log\left(
\frac{V'(\theta) V^\alpha(\theta) \prod_{m =
N_F+1}^{n-1}|\sqrt{V(\theta)} + a_m|
\theta_\varepsilon^\alpha}
{V'(\theta_\varepsilon) V^\alpha(\theta_\varepsilon) 
\prod_{m = N_F+1}^{n-1}|\sqrt{V(\theta_\varepsilon)} +
a_m|
\theta^\alpha}\right) < (-1 + 2\varepsilon)(s - s_o)\ ;$$
therefore
$$\frac{V'(\theta)}{V'(\theta_\varepsilon)} < \frac{\hat M_{\varepsilon}}{
V(\theta)^\alpha \prod_{m =
N_F+1}^{n-1}|\sqrt{V(\theta)} + a_m|}
\frac{e^{-(1-2\varepsilon) s}}
{e^{-(1-2\varepsilon) s_o}} = \tilde M_{\varepsilon} \left(\frac{1 -
\theta}{1-\theta_\varepsilon}\right)^{1-2\varepsilon}$$
for some suitable positive constant $\tilde M_{\varepsilon}$. 
Since for any $s$,  $V'(\theta(s)) = \dot V(\theta) (1-\theta)$, it follows that
$$
\frac{\dot V(\theta) (1-\theta) }{\dot
V(\theta_\varepsilon)(1-\theta_\varepsilon)} < \tilde
M_{\varepsilon} 
\left(\frac{1 - \theta}{1-\theta_\varepsilon}\right)^{1 - 2\varepsilon}$$
and  this implies the claim.\par
\bigskip
\bigskip
We have now all ingredients to complete the proof of Theorem 5.5. In fact, 
 consider the solution $V$ obtained by Claim 3. If we can show that 
$$\lim_{\theta\to 0} \dot V(\theta) = \dot V_o >0\tag 5.23$$
and that   
$\lim_{\theta \to 1} \dot V(\theta)$, $\lim_{\theta \to 1}
\overset{\cdot \cdot}\to V (\theta)$
and $\lim_{\theta \to 1} \overset{\cdot \cdot \cdot}\to V(\theta)$ 
exist and are finite, we are done.\par
First of all,  it is not hard to check that, by   (5.10), 
$\log\left(\dot V(\theta) \left(\frac{V(\theta)}{\theta}\right)^\alpha\right)$ 
verifies the Cauchy condition for $\theta$ tending to $0$. 
Hence 
$\lim_{\theta \to 0} \log\left(\dot V(\theta)
\left(\frac{V(\theta)}{\theta}\right)^\alpha\right)$  
exists and, using again (5.10),  one can check that this limit is positive. 
 Since $\lim_{\theta \to 0} V(\theta) = 0$,  by 
de L'H\^opital Theorem,  we also have that
$\lim_{\theta \to 0} \left(\frac{V(\theta)}{\theta}\right)^{\alpha+1}$ exists and
that 
it is positive. 
From these facts,  we conclude that  (5.23) is verified.\par
\bigskip
Consider now a value $\varepsilon < \frac{1}{2}$ and let $M_{\varepsilon}$ be the
constant 
given in Claim 5. In this case, 
$$0 < \lim_{\theta \to 1^-} \dot V (1-\theta)^{\frac{1}{2} + \varepsilon} 
= \left[\lim_{\theta \to 1^-}  \dot V(1 - \theta)\right] (1-\theta)^{-
\frac{1}{2} +\varepsilon} < 
 M_\varepsilon
\lim_{\theta \to 1^-} (1- \theta)^{\frac{1}{2} - \varepsilon} = 0\ .\tag 5.24$$
From (5.24) and de L'H\^opital theorem, we get  
$$\lim_{\theta\to 1^-} \frac{2(1 + \alpha') - \hat c \sqrt{\frac{V}{\theta}}}{(1
-
\theta)^{\frac{1}{2}-\varepsilon}} = 
\lim_{\theta\to 1^-} \frac{2(1 + \alpha')\sqrt{\theta} - \hat c \sqrt{V}} {
(1 -
\theta)^{\frac{1}{2}-\varepsilon}} = $$
$$ = 
\lim_{\theta\to 1^-} \frac{\frac{2(1 + \alpha')}{2 \sqrt{\theta}} - \frac{\hat c
\dot
V}{2 \sqrt{V}}}
{(\frac{1}{2}+\varepsilon)(1 - \theta)^{-(\frac{1}{2}+\varepsilon)}} = $$
$$ = 
 \lim_{\theta\to 1^-} \frac{2(1 + \alpha')}{2 (\frac{1}{2}+\varepsilon)
\sqrt{\theta}}(1-
\theta)^{\frac{1}{2}+\varepsilon}
- \frac{\hat c^2}{4(1 +\alpha') (\frac{1}{2}+\varepsilon)} \lim_{\theta\to 1^-}
\dot V (1 -
\theta)^{\frac{1}{2}+\varepsilon} = 0\ .\tag 5.25$$
Now, using (5.24), (5.25) and (5.1), we obtain that for all $\theta > \theta_o$,
for some suitable 
$\theta_o$, 
there exists a positive constant 
$C_1>0$ such that
$$\left|\overset{\cdot \cdot}\to V (1- \theta)\right| = 
 \left|- \frac{ \dot V^2 (1 - \theta) }{2\sqrt{V}}
\left(\frac{ 2 \alpha }{\sqrt{V}} + 
\sum_{m = N_F+1}^{n-1} \frac{1}
{\sqrt{V} + a_m}\right) +\right. $$
$$ + \left.\left( 
2(1 + 
\alpha')  + \frac{\alpha}{\theta} - \hat c\sqrt{\frac{V}{\theta}}\right) \dot
V\right| <
\dot V \left( C_1 (1 - \theta)^{\frac{1}{2} - \varepsilon} +
\frac{\alpha}{\theta}\right)\tag 5.26$$
Therefore, by integration, we get that for any two $\theta_o <\theta_1 <\theta_2
< 1$
$$\log\left(\frac{\dot V_{\theta_2}}{\dot V_{\theta_1}}\right) < C_1
\int_{\theta_1}^{1} \frac{1}{(1 - \theta)^{\frac{1}{2} + \varepsilon}} d\theta 
+ \int_{\theta_1}^{1}\frac{\alpha}{\theta} d\theta = $$
$$ = 
\frac{2 C_1}{ 1 - 2 \varepsilon}(1 - \theta_1)^{\frac{1}{2} - \varepsilon} -
\log(\theta_1^\alpha) = 
C_2 (1 - \theta_1)^{\frac{1}{2} - \varepsilon} - \log(\theta_1^\alpha) \ .\tag
5.27$$
and
$$\log\left(\frac{\dot V_{\theta_2}}{\dot V_{\theta_1}}\right) > - C_2 (1 -
\theta_1)^{\frac{1}{2} - \varepsilon} +
\log(\theta_1^\alpha)\ .\tag 5.28$$
From (5.27) and (5.28), it  follows  that $\log(\dot V_{\theta})$ 
verifies the Cauchy condition for $\theta$ tending to $1$ and hence 
that $\lim_{\theta \to 1^-}\dot V_{\theta}$
exists  and it is
finite. Moreover by (5.10) this limit is positive as required.\par
The claims that $\lim_{\theta \to 1}
\overset{\cdot \cdot}\to V (\theta)$
and $\lim_{\theta \to 1} \overset{\cdot \cdot \cdot}\to V(\theta)$ 
exist and are finite can be checked immediately from (5.1), with a
straightforward 
application of  de L'H\^opital theorem. \qed
\enddemo
\bigskip
\demo{Proof of Theorem 5.1} One can check that 
for any K-manifold described in Table 1 of Corollary
3.5, with the only exceptions 
of the manifold in n.1 with $G = \SU_3$ and $F = \CP^2$ and 
those in n.2 
with $G = \SU_p \times \SU_2$, $p>2$, the condition
d) of Proposition 5.2 is satisfied. Then
the conclusion follows as 
direct corollary of Proposition 5.2 and Theorem 5.6.
\qed
\enddemo
\bigskip
\bigskip
%\newpage
\head APPENDIX
\endhead
\bigskip
In the following table, we adopt the same notation for simple root system
adopted in \cite{GOV} (see also \cite{AS}). \par
For any case of Table 1, we list the Lie algebra $\g$, the Morimoto-Nagano
algebra 
$\g_F$, the value of $N_F$ (which in all considered cases 
coincides with $N_F^{(1)}$, being $N_F^{(2)} = 0$), the 
element $\theta_D$ in the dual space $\h^*$ of a Cartan subalgebra of 
$\g$,  which is  $\B$-dual to the element $- iZ_\D$, the set of roots in 
$R'_+$, the element $\theta^\kappa \in \h^*$, which is 
$\B$-dual to $ - i Z^\kappa$ and all values for 
$$c a_{\beta_m} = \frac{\B(Z^\kappa_\perp, i H_\beta)}
{\B(Z_\D, i H_\beta)}$$
which occurs for $\beta_m \in R'$.
Notice that for all considered cases
$Z^\kappa$ is orthogonal to $Z_\D$ and hence that 
 $Z^\kappa_\perp = Z^\kappa$.
\bigskip 
\moveleft 0.3 cm
\vbox{\offinterlineskip
\halign {\strut
\vrule\hfil\ $#$\ \hfil
 &
\vrule\hfil\ $#$\ \hfil
&
\vrule\hfil\ $#$\ \hfil
&
\vrule\hfil\ $#$\ \hfil
&
\vrule\hfil\ $#$\ \hfil
&
\vrule\hfil\ $#$\ \hfil
&
\vrule\hfil\ $#$\ \hfil
&
\vrule\hfil\ $#$ \hfil \vrule
\cr \noalign{\hrule} 
n^o 
&
\phantom{\frac{\frac{1}{1}}{\frac{1}{1}}}\g
\ \ 
&
\g_F
&
N_F
&
\theta_\D
&
R'_+
&
\theta^\kappa
&
\underset{\phantom{B}}\to{
\overset{\phantom{A}}\to{
\frac{\B(Z^\kappa_\perp, i H_\beta)}
{\B(Z_\D, i H_\beta)}
}}
\cr \noalign{\hrule}
1 
&
\underset{\phantom{B}}\to{
\overset{\phantom{A}}\to{
\su_{\ell+1}
}}
&
\su_2
&
1
&
 - \frac{1}{2}(\varepsilon_1 - \varepsilon_2)
&
\underset{\phantom{B}}\to{
\overset{\phantom{A}}\to{
\smallmatrix
\varepsilon_1 - \varepsilon_a,
\ \varepsilon_2 - \varepsilon_a\\
\phantom{a}\\
 3\leq a\leq \ell + 1\\
\phantom{a}
\endsmallmatrix
}}
&
\smallmatrix 
(\ell - 1)(\varepsilon_1 + \varepsilon_2) -
\\
\phantom{a} 
\\
- 2\sum_{a = 3}^{\ell+1}\varepsilon_a 
\endmatrix
&
\smallmatrix
\pm 2(\ell + 1)
\endsmallmatrix
\cr \noalign{\hrule}
2 
&
\underset{\phantom{B}}\to{
\overset{\phantom{A}}\to{
\matrix
\su_{p+1} \oplus\\
 \su_{q+1}
\\
\ 
\\
\smallmatrix
(p + q > 2)
\endsmallmatrix
\endmatrix
}}
&
\smallmatrix 
\su_2 \oplus \su_2
\endsmallmatrix
&
2
&
\matrix
 - \frac{1}{2}(\varepsilon_1 - \varepsilon_2) - \\
- \frac{1}{2}(\varepsilon'_1 - \varepsilon'_2)
\endmatrix
&
\underset{\phantom{B}}\to{
\overset{\phantom{A}}\to{
\matrix 
\smallmatrix
\varepsilon_1 - \varepsilon_a,
\ \varepsilon_2 - \varepsilon_a\\
\phantom{a}\\
 3\leq a\leq p + 1\\
\phantom{a}
\endsmallmatrix
\\
\smallmatrix
\varepsilon'_1 - \varepsilon'_b,
\ \varepsilon'_2 - \varepsilon'_b\\
\phantom{a}\\
 3\leq b\leq q + 1\\
\phantom{a}
\endsmallmatrix
\endmatrix 
}}
&
\smallmatrix 
(p - 1)(\varepsilon_1 + \varepsilon_2) 
+ \\
\phantom{a}
\\
 + (q - 1)(\varepsilon'_1 + \varepsilon'_2) - \\
\phantom{a}
\\
- 2\sum_{a = 3}^{p+1}\varepsilon_a - \\
\phantom{a}
\\
- 2\sum_{b = 3}^{q+1}\varepsilon'_b 
\endsmallmatrix
&
\smallmatrix
\pm 2(p+1)\ ,\\ 
\pm 2(q + 1)
\endsmallmatrix 
\cr \noalign{\hrule}
3
&
\underset{\phantom{B}}\to{
\overset{\phantom{A}}\to{
\matrix
\su_{\ell+1}\\
\ 
\smallmatrix
\ell > 3
\endsmallmatrix
\endmatrix
}}
&
\so_6
&
4
&
\matrix
 - \frac{1}{2}(\varepsilon_1 + \varepsilon_2 - \\
- \varepsilon_3 - \varepsilon_4)
\endmatrix 
&
\underset{\phantom{B}}\to{
\overset{\phantom{A}}\to{
\smallmatrix
\varepsilon_i - \varepsilon_a\\
\phantom{a}\\
1 \leq i \leq 4 \\
\phantom{a}\\
 5\leq a\leq \ell + 1\\
\phantom{a}
\endsmallmatrix
}}
&
\smallmatrix 
(\ell - 1)\sum_{i = 1}^4\varepsilon_i -
\\
\phantom{a}
\\
- 4\sum_{a = 5}^{\ell+1}\varepsilon_a 
\endsmallmatrix
&
\smallmatrix
\pm (2 (\ell - 1) + 8)
\endsmallmatrix
\cr \noalign{\hrule}
4
&
\underset{\phantom{B}}\to{
\overset{\phantom{A}}\to{
\so_{10}
}}
&
\so_8
&
6
&
 \frac{1}{2}\sum_{i = 2}^5
\varepsilon_i
&
\underset{\phantom{B}}\to{
\overset{\phantom{A}}\to{
\smallmatrix
\varepsilon_1 \pm \varepsilon_i\\
\phantom{a}\\
2\leq i,j\leq 5\\
\phantom{a}
\endsmallmatrix
}}
&
\smallmatrix 
8 \varepsilon_1
\endsmallmatrix
&
\smallmatrix
\pm 16
\endsmallmatrix
\cr \noalign{\hrule}
5 
&
\underset{\phantom{B}}\to{
\overset{\phantom{A}}\to{
\goth e_6
}}
&
\so_{10}
&
8
&
\underset{\phantom{B}}\to{
\overset{\phantom{A}}\to{
\matrix
 - \frac{1}{2}(2\varepsilon_1 + \varepsilon_6\\
\phantom{a}\\
+ \varepsilon)
\endmatrix
}}
&
\underset{\phantom{B}}\to{
\overset{\phantom{A}}\to{
\smallmatrix
\varepsilon_i - \varepsilon_6 , \\
\phantom{a}\\
\varepsilon_i +\varepsilon_j + \varepsilon_k + 
\varepsilon , \\
\phantom{a}\\
2 \varepsilon\\
\phantom{a}\\
1\leq i,j,k\leq 5\\
\phantom{a}
\endsmallmatrix
}}
&
\smallmatrix 
12(- \varepsilon_6 + \varepsilon)
\endsmallmatrix
&
\smallmatrix
\pm 24
\endsmallmatrix
\cr \noalign{\hrule}
}}
\centerline{\bf Table A1}
\bigskip
\bigskip
We now list the explicit expressions for the 
integrals $\int_0^{V_{\hat c}}g(x) dx$, which appear in the proof of Lemma 5.7,
for all cases of Table A1. We normalize $V_{\hat c} = 1$, i.e. we assume
that  $\hat c = 2(1+\alpha')$ 
and $c = -{{\varepsilon_F \hat c}\over {\B(Z_\D,Z_\D)}}$. \par
If $\tilde M$ is as in  n.1 with $\varepsilon_F = 1$, we have 
$${1\over {4^{l-1}}} \int_0^1 ((l+1)^2 - 4x)^{l-1} (1 - 2\sqrt{x}) dx\ ,$$
which is negative, as it is proved in [GC]. With $\varepsilon_F = 2$, we obtain 
$${1\over {2 \cdot 9^{l-1}}} \int_0^1 ((l+1)^2 - 9x)^{l-1} (1 - 3\sqrt{x}) dx\ .$$
Now, 
$$\eqalign{ \int_0^1 &((l+1)^2 - 9x)^{l-1}   (1 - 3\sqrt{x}) dx < \cr
{} & < (l+1)^{2(l-1)} \left(\int_0^{1/9} 1 - 3\sqrt{x} dx\right) + ((l+1)^2 - 9)^{l-1}
\left(\int_{1/9}^1 1 - 3\sqrt{x} dx\right) = \cr
{}& = {1\over{27}}((l+1)^{2(l-1)} - 28 ((l+1)^2 - 9)^{l-1})\ ,}$$
which is negative for all $l \geq 3$, by the same arguments in  [GC].\par
If $\tilde M$ is as in  n.2,  the only possibility for $\varepsilon_F$ is  
$\varepsilon_F = 2$ and the integral is 
$${1\over{4^{p-1}4^{q-1}}} \int_0^1 \sqrt{x} ((p+1)^2 - 4x)^{p-1} 
((q+1)^2 - 4x)^{q-1} (1 - 2\sqrt{x}) dx\ .$$
Now, 
$$\eqalign{
\int_0^1 & \sqrt{x} ((p+1)^2 - 4x)^{p-1} 
((q+1)^2 - 4x)^{q-1} (1 - 2\sqrt{x}) dx < \cr
{} & < (p+1)^{2(p-1)}(q+1)^{2(q-1)}\left(\int_0^{1/4} \sqrt{x} (1 - 2\sqrt{x})dx\right) + \cr
{}& + 
((p+1)^2 - 4)^{p-1}((q+1)^2 - 4)^{q-1}\left(\int_{1/4}^1 
\sqrt{x} (1 - 2\sqrt{x})dx\right) =\cr
{}& = {1\over {48}}\left( (p+1)^{2(p-1)}(q+1)^{2(q-1)} - 17 
((p+1)^2 - 4)^{p-1}((q+1)^2 - 4)^{q-1}\right)\ ,}$$
which is again negative for all $p,q \geq 2$.\par
If $\tilde M$ is as in  n.3, we have only the possibility $\varepsilon_F = 2$ and the integral is
$${1\over{9^{2l-6}}} \int_0^1 x^{3/2} ((l+3)^2 - 9x)^{2(l-3)} (2 - 3\sqrt{x}) dx,$$
which can be proved to be negative for all $l\geq 4$.\par
If $\tilde M$ is as in  case n.4, with $\varepsilon_F = 2$ we have that the integral is  
$$\int_0^1 x^{5/2} (4 - x)^4 (3 - 4\sqrt{x}) dx < 0,$$
as it can be  verified using e.g. some symbolic manipulation computer program like
{\it Mathematica\/}. If we assume that $\varepsilon_F = 1$, the integral is 
$${1\over{7^8}} \int_0^1 x^{5/2}(256 - 49x)^4 (6 - 7\sqrt{x}) dx\ ,$$
which can be checked to be positive.\par
If $\tilde M$ is as in  n.5, with $\varepsilon_F = 1$, the integral is 
$$\int_0^1 x^{7/2} ({{64}\over 9} - x)^8 (8 - 9\sqrt{x}) dx\ ,$$
which can be checked to be positive;
in case $\varepsilon_F = 2$, the integral is 
$$\int_0^1 x^{7/2} ({{144}\over{25}} - x)^8 (4 - 5\sqrt{x}) dx$$
which can be checked to be negative. 
\par
  
\bigskip

\hbox{\parindent=0pt\parskip=0pt
\vbox{\hsize=2.7truein
\obeylines
{
Fabio Podest\`a
Dip.~di Matematica e Appl.~per l'Arch.
Universit\`a di Firenze
P.zza Ghiberti, 27
I-50142 Firenze
Italy
}\medskip
podesta\@math.unifi.it
}\hskip 1.5truecm
\vbox{\hsize=3.7truein
\obeylines
{
Andrea Spiro
Dip. di Matematica e Fisica
Universit\`a di Camerino
Via Madonna delle Carceri
I-62032 Camerino (MC)
Italy
}\medskip
spiro\@campus.unicam.it
}
}

\enddocument

\bigskip
\bigskip



\Refs  
\widestnumber\key{GOV}  
  

\ref
\key Ah
\by D. N. Ahiezer
\paper Algebraic groups acting transitively
in the complement of a homogeneous 
hypersurface
\jour Soviet Math. Dokl.
\vol 20
\yr 1979
\pages 278--281
\endref

\ref
\key Al
\by D. V. Alekseevsky
\paper Flag Manifolds
\inbook Sbornik Radova, 11 Jugoslav. Seminr.
\vol 6
\issue 14
\yr 1997
\publ Beograd
\pages 3--35
\endref

\ref  
\key AA  
\by A.V. Alekseevsky, D.V.Alekseevsky  
\paper G-manifolds with one dimensional orbit space  
\jour Adv. in Sov. Math.  
\vol 8  
\yr 1992   
\pages 1--31  
\endref 

\ref  
\key AA1  
\bysame   
\paper Riemannian G-manifolds with one dimensional orbit space  
\jour Ann. Glob. Anal. and Geom.  
\vol 11  
\yr 1993  
\pages 197--211  
\endref  

\ref
\key AP
\by D. V. Alekseevsky and A. M. Perelomov
\paper Invariant Kahler-Einstein metrics
on compact homogeneous spaces
\jour Funktsional. Anal. i Prilozhen.
\vol 20
\yr 1986 
\issue 3
\transl\nofrills Engl. transl. in
\jour Funct. Anal. Appl.  
\vol 20
\yr 1986
\issue 3
\pages 171--182
\endref

\ref
\key AS
\by D.V. Alekseevsky, A. Spiro
\paper Invariant CR-structures on compact homogeneous manifolds
\paperinfo Preprint at the Los Alamos Electronic Archive 
(http://xxx.lanl.gov/math.DG/9904054)
\yr 1999
\endref

\ref 
\key Ber
\by L. B\'erard Bergery
\paper Sur de nouvelles vari\'et\'es riemanniennes d'Einstein
\jour Instute Elie Cartan
\vol 6
\yr 1983
\pages 1--60
\endref

\ref
\key Be 
\by A.L. Besse
\book Einstein manifolds
\publ Springer-Verlag 
\yr 1986
\endref

\ref
\key BFR
\by M. Bordermann, M. Forger and H. R\"omer
\paper Homogeneous K\"ahler Manifolds: paving the way
towards new supersymmetric Sigma Models
\jour Comm. Math. Phys.
\vol 102
\yr 1986
\pages 605--647
\endref
 
\ref  
\key Br  
\by G.E. Bredon  
\book Introduction to compact transformation groups  
\publ Acad. Press N.Y. London  
\yr 1972  
\endref



\ref
\key DW
\by A. Dancer, M. Wang
\paper K\"ahler-Einstein metrics of cohomogeneity one
\publ Math. Ann. 
\yr 1998
\vol 312
\pages 503--526
\endref

\ref 
\key GC
\by D. Guan, X. Chen 
\paper Existence of Extremal Metrics on Almost
Homogeneous Manifolds of Cohomogeneity One
\jour Asian J. Math. 
\vol 4
\yr 2000
\pages 817--830
\endref


\ref
\key GOV
\by V. V. Gorbatsevic, A. L. Onishchik and E. B. Vinberg
\paper Structure of Lie Groups and Lie Algebras
\inbook in Encyclopoedia of Mathematical Sciences -
Lie Groups and Lie Algebras III
\ed A. L. Onishchik and E. B. Vinberg
\publ Springer-Verlag -- VINITI
\publaddr Berlin
\yr 1993 (Russian edition: VINITI, Moscow,1990) 
\endref

\ref 
\key HO
\by A. Huckleberry, E. Oeljeklaus
\paper A characterization of complex homogeneous cones
\jour Math. Z.
\vol 170
\yr 1980
\pages 181--194
\endref

\ref
\key HS
\by A. Huckleberry, D. Snow
\paper Almost-homogeneous K\"ahler manifolds with hypersurface 
orbits
\jour Osaka J. Math.
\vol 19
\yr 1982
\pages 763--786
\endref

\ref
\key KS
\by N. Koiso, Y. Sakane
\paper Non-homogeneous K\"ahler-Einstein metrics on compact complex manifolds II
\jour Osaka J. Math.
\vol 25 
\yr 1988
\pages 933--959
\endref 

\ref
\key KW
\by J. L. Kazdan, F. W. Warner
\paper Curvature functions for open 2-manifolds
\jour Ann. of Math. 
\vol 99
\yr 1974
\pages 203--219
\endref 

\ref
\key MN
\by A. Morimoto, T. Nagano
\paper On pseudo-conformal transformations of hypersurfaces
\jour J. Math. Soc. Japan
\vol 15 
\yr 1963
\pages 289--300
\endref 

\ref
\key Oe
\by E. Oeljeklaus
\paper Fast homogene K\"ahlermannigfaltigkeiten
mit verschwindender erster Bettizahl
\jour manuscripta math.
\vol 7
\yr 1972
\pages 175--183
\endref

\ref
\key Pa
\by R. S. Palais
\paper On the existence of slices for actions of non-compact Lie groups
\jour Ann. of Math.
\issue 2
\vol 73
\yr 1961
\pages 480--484
\endref


\ref
\key PS
\by F.Podest\`a, A. Spiro
\paper K\"ahler manifolds with large isometry group
\jour Osaka J. Math.
\vol 36
\yr 1999
\pages 805--833
\endref

\ref
\key Sa
\by Y. Sakane
\paper Examples of compact K\"ahler Einstein manifolds with positive 
Ricci curvature 
\jour Osaka J. Math.
\vol 23
\yr 1986
\pages 585--616
\endref

\ref
\key Sp
\by A. Spiro
\paper The Ricci tensor of an almost homogeneous K\"ahler manifold
\paperinfo Preprint at the Los Alamos Electronic Archive 
(http://xxx.lanl.gov/math.DG/0101172)
\yr 2001
\endref

\ref
\key Uc
\by F. Uchida
\paper Classification of compact transformation groups
on cohomology complex projective spaces with 
codimension one orbits
\jour Japan. J. Math.
\vol 3
\yr 1977
\pages 141--189
\endref

\ref
\key Ve
\by L. Verdiani
\paper Invariant metrics on cohomogeneity one manifolds
\jour Geom. Dedicata
\vol 77
\yr 1999
\pages 77--111
\endref
\endRefs
\bye